\documentclass[%
 aip,
 amsmath,amssymb,
 reprint,%
]{revtex4-1}

\usepackage{graphicx}
\usepackage{dcolumn}
\usepackage{bm}

\usepackage[utf8]{inputenc}
\usepackage{algorithm}
\usepackage[T1]{fontenc}
\usepackage{mathptmx}
\usepackage{etoolbox}
\usepackage{braket}
\usepackage{algpseudocode}

\makeatletter
\def\@email#1#2{%
 \endgroup
 \patchcmd{\titleblock@produce}
  {\frontmatter@RRAPformat}
  {\frontmatter@RRAPformat{\produce@RRAP{*#1\href{mailto:#2}{#2}}}\frontmatter@RRAPformat}
  {}{}
}%
\makeatother
\begin{document}

\preprint{AIP/123-QED}

\title[Reduced basis algorithm for solving nonlinear differential equations on quantum computers]{Reduced basis algorithm for solving nonlinear differential equations on quantum computers}
\author{Monica L\u{a}c\u{a}tu\c{s}}
\email{m.i.l.lacatus@tudelft.nl}
\affiliation{Delft Institute of Applied Mathematics, Delft University of Technology, Mekelweg 4, 2628 CD Delft, The Netherlands}

\author{Matthias M\"{o}ller}
\affiliation{Delft Institute of Applied Mathematics, Delft University of Technology, Mekelweg 4, 2628 CD Delft, The Netherlands
}

\author{Sauro Succi}%
\affiliation{ Fondazione Istituto Italiano di Tecnologia, Center for Life Nano-Neuroscience
at la Sapienza, Viale Regina Elena 291, 00161 Roma, Italy
}%


\begin{abstract}
As quantum computing moves toward scientific computing applications, nonlinear differential equations remain a central challenge since quantum evolution is intrinsically linear. Existing quantum algorithms typically handle nonlinear updates using quantum state copies, mean-field constructions, or linearization techniques, which can introduce exponential copy growth, copy-number-dependent approximation errors, or restrictive convergence and stability conditions. In this work, we introduce a reduced basis algorithm (RBA) for polynomial nonlinear ordinary differential equations (ODEs) and spatially discretized partial differential equations (PDEs). After time discretization, the method composes the resulting polynomial update map over $m$ timesteps, identifies the reduced monomial basis appearing in this composed map, and constructs a linear RBA operator whose action recovers the exact $m$-timestep nonlinear dynamics. Thus, at the level of the chosen discrete update rule, the method introduces no additional approximation error beyond the time discretization error. The qubit number requirement is governed by the size of the reduced monomial basis. For an $n$-dimensional polynomial ODE system of degree $p>1$, the lifted register requires at most $q_m^{\mathrm{ODE}} = O(nm\log p)$ qubits in the full basis scenario. For PDEs discretized on $N^D$ grid points, a locality-based construction requires at most $q_m^{\mathrm{PDE}} = O(D\log N + n m^{D+1}\log p)$ qubits. Hence, the dependence on the grid size remains logarithmic, while the nonlinear overhead is controlled by local reduced basis size. The main computational burden is moved from the quantum computer to a classical preprocessing step, where the reduced monomial basis and RBA operator are constructed for the chosen timestep window. Through numerical tests on the Lorenz system and the one-dimensional Burgers equation, we verify that the RBA reproduces the corresponding discrete time nonlinear dynamics exactly, while exposing the trade-off between timestep composition, reduced basis growth, and locality.
\end{abstract}

\maketitle

\section{Introduction}
\label{sec:intro}

Nonlinear differential equations arise throughout science and engineering, from fluid dynamics and plasma physics to nonlinear optics and reaction-diffusion systems. Motivated by the prospect of computational speedups and memory reductions, there is growing interest in developing efficient quantum algorithms for such systems. However, nonlinearity remains one of the central obstacles, raising fundamental questions about whether nonlinear dynamics can be efficiently represented within the intrinsically linear framework of quantum computation.

One of the earliest quantum algorithms for nonlinear differential equations was proposed by Leyton and Osborne \cite{leyton2008quantum}. Their method targets sparse systems of ODEs whose nonlinear terms are polynomials. The continuous time problem is first discretized using the explicit Euler method and the solution vector is encoded in the amplitudes of a quantum state. Since quantum amplitudes cannot be transformed nonlinearly by a unitary operation, the algorithm uses multiple copies of the state to represent the polynomial terms appearing in the update. For example, if two copies of a quantum state $\lvert \phi\rangle\lvert \phi\rangle$ are available, their tensor-product
state contains amplitudes proportional to products of the form $z_k z_l$. These amplitudes can then be used to represent the quadratic monomials needed to evaluate quadratic polynomial nonlinearities. Cubic nonlinearities would similarly require three copies of the quantum state, and so on. The algorithm then constructs a linear operator on the enlarged tensor-product space and embeds it into a Hamiltonian simulation. This gives a nondeterministic quantum implementation of an explicit Euler step. However, each nonlinear transformation consumes multiple copies of the current state and succeeds only probabilistically. The number of qubits required for this algorithm scales polylogarithmically with the number of variables, but exponentially with the number of timesteps. This makes the method impractical for long-time simulations.

Several years later, Lloyd et al. \cite{lloyd2020quantum} proposed a quantum algorithm for nonlinear differential equations that avoids the exponential copy overhead of the approach of Leyton and Osborne \cite{leyton2008quantum}. Polynomial nonlinearities are evaluated by introducing $N_c$ identical copies of the state and applying weak, symmetrized interactions among them. In the mean-field limit, the reduced state of any single copy evolves approximately according to the desired nonlinear equation. Thus, unlike the method of Leyton and Osborne \cite{leyton2008quantum}, the algorithm does not consume fresh copies at each time step. For general non-Hermitian nonlinear differential equations, Lloyd et al. \cite{lloyd2020quantum} combine this mean-field construction with a quantum linear solver algorithm (QLSA). After time discretization, the nonlinear evolution is embedded into a linear system whose solution is a quantum history state containing the approximate solution at all time steps. Under suitable stability and approximation assumptions, the qubit number scaling can be polynomial in the evolution time. In the first-order Euler/Trotter analysis, the mean-field error is suppressed as $O(1/N_c)$ but accumulates over time. The main limitation is that the approximation can fail when the nonlinear dynamics strongly amplify perturbations, for example in systems with large positive Lyapunov exponents. Tennie and Palmer \cite{tennie2023quantum} illustrate this accumulation numerically for the toy model $\dot{x}=x-\alpha x^3$, varying both the initial condition and nonlinear coupling strength $\alpha$.

A separate linearization-based approach avoids the use of multiple quantum state copies by embedding the nonlinear dynamics into an infinite-dimensional linear system through Carleman linearization. In this approach, the original state $x(t)$ is lifted to a hierarchy of monomials, represented by tensor powers such as $x(t),x(t)^{\otimes 2}, x(t)^{\otimes 3}$, and so forth. After truncation at a finite Carleman order, the resulting system is linear, and the degree-one component of the lifted state—the block corresponding to $x(t)$-approximates the solution of the original nonlinear system. This makes the method attractive for quantum algorithms, since the truncated linear system can be treated using QLSA and amplitude encoding can retain logarithmic dependence on the original state dimension. However, the lifted dimension grows rapidly with the truncation order, so efficiency depends on keeping this order small.

Liu et al. \cite{liu2021efficient} proposed one of the first rigorous quantum algorithms based on Carleman linearization for dissipative quadratic nonlinear ODEs. Their main assumption is that the nonlinear and forcing terms are sufficiently weak compared with the linear dissipation. Under this condition, the Carleman truncation converges, and the resulting finite-dimensional linear system can be solved using QLSA based on an Euler time discretization. This gives polylogarithmic scaling in the original dimension and polynomial scaling in the evolution time and inverse error. However, this efficiency is conditional. If the dissipation is too weak, or if the nonlinear terms dominate, the advantage can disappear. 

Krovi \cite{krovi2023improved} improved this framework by replacing the Euler-based linear ODE solver with a higher-order quantum linear differential equation algorithm. This gives an exponential improvement in the dependence on the target error, even when the Carleman linearized matrix is not diagonalizable, and allows a broader class of linear dissipative parts with negative log-norm to be simulated.

Costa et al. \cite{costa2025further} further developed the quantum Carleman approach, combining higher-order time discretization with a tighter analysis of the Carleman truncation error. Their work also highlights an important measurement issue which is that the Carleman state contains many lifted components, while the physically relevant solution is contained only in the first block corresponding to $x(t)$. Without rescaling, the higher-order lifted components can dominate the normalized Carleman state, so the probability of extracting the physical block can become exponentially small as the Carleman truncation order increases. This is important because the truncation order must be increased to obtain a more accurate approximation of the nonlinear dynamics. For discretized PDEs, the problem can be worse because the size of the solution vector also grows with the number of spatial degrees of freedom. Rescaling can improve the success probability, but it also changes the balance between nonlinear growth and linear dissipation in the lifted system. Thus, even if the original nonlinear equation is stable, the rescaled Carleman system may fail to satisfy the stability conditions needed for an efficient quantum algorithm.

Carleman linearization has been used to construct quantum algorithms for several methods and models in fluid dynamics, including the lattice Boltzmann method \cite{Sanavio2024,Jennings2025} and reaction-diffusion equations \cite{Liu2023ReactionDiffusion,An2026QuantumAlgorithms}. These works extend the same basic idea to more concrete applications, but they do not overcome a central limitation of Carleman linearization which is that finite-order Carleman truncations generally provide accurate approximations only within a restricted convergence regime. When the dynamics remain within this regime, the truncated lifted system can approximate the original nonlinear evolution well. However, once the trajectory leaves this regime the approximation can fail. Novikau and Joseph \cite{novikau2025globalizing} emphasize this point and propose globalized, piecewise Carleman methods that switch between local linearization charts. These methods can handle systems that are inaccessible to standard Carleman linearization, including systems with multiple equilibria, limit cycles, and chaotic attractors. However, it is not yet clear whether and how this procedure can be implemented as a quantum algorithm.

Overall, the efficiency of Carleman-based quantum algorithms is best understood as regime dependent, with the most favorable behavior arising for weakly nonlinear, stable systems whose trajectories remain in a locally convergent region. Other linearization-based approaches, including Koopman linearization \cite{katz2025efficient}, Fokker-Planck formulations \cite{tennie2024solving}, and quantum homotopy perturbation methods \cite{xue2021quantum,xue2022quantum}, have also been analyzed, but they generally inherit similar limitations to Carleman, including truncation overhead, conditioning issues, and restrictive convergence or stability requirements.

Variational quantum algorithms (VQAs) have also been investigated as an alternative approach for solving nonlinear differential equations. Notable examples include the work of Lubasch et al. \cite{lubasch2020variational}, who encode the solution of the nonlinear differential equation as a variational quantum state and estimate nonlinear contributions by measuring observables on multiple copies of the state. Kyriienko et al. \cite{kyriienko2021solving}, in contrast, approximate the solution using parametrized quantum circuits and formulate the problem as the minimization of a differential-equation residual, with derivatives evaluated analytically from the circuit structure. These approaches are attractive because they can exploit quantum state representations while remaining compatible with hybrid optimization. However, their overall complexity remains difficult to assess, since general analytic convergence guarantees for the variational optimization are not yet established. 

Recent work has also connected VQAs with tensor-network representations. Siegl et al. \cite{siegl2026tensorprogrammable} propose a tensor-programmable quantum circuit framework in which discretized differential operators are represented as matrix product operators and mapped to quantum circuits. This gives a flexible tensor-network-based formulation of variational quantum time stepping. Its efficiency, however, still depends on ansatz expressivity, optimization performance, and the sampling cost associated with post-selection and cost-function estimation.

In this work, we introduce a reduced basis algorithm (RBA) for solving polynomial nonlinear differential equations on quantum computers. This algorithm inherits the monomial basis construction of Carleman linearization, but it does not construct an infinite-dimensional system that requires truncation at a specific order. Instead, after discretizing the continuous time problem, we represent the $m$-timestep nonlinear evolution as a composed polynomial map acting on the initial condition. The nonlinear update is then evaluated by applying a linear operator to a reduced basis containing only the monomials that actually appear in this composed map. In this way, the algorithm recovers the exact fully discrete nonlinear time evolution, with no additional approximation error beyond the error introduced by the chosen time discretization scheme.

In this algorithm, the main computational overhead is shifted from the quantum computer to a classical pre-processing stage. For a fixed timestep window, the discrete nonlinear update map is composed classically, the reduced monomial basis appearing in the composed map is identified, and the associated linear operator is constructed. This precomputed linear operator is then block-encoded and applied to a lifted quantum state containing the monomials in the reduced basis. Therefore, the principal cost is governed by the size of the reduced monomial basis. For spatially discretized PDEs, this cost can be reduced by exploiting stencil locality, such that the reduced lifted basis and linear operators are constructed locally rather than over the full global state.

The remainder of the paper is organized as follows. Section \ref{sec:RBA} introduces the RBA construction for systems of polynomial nonlinear ordinary differential equations (ODEs) and extends the same idea to polynomial nonlinear partial differential equations (PDEs). Section \ref{sec:QAC} discusses quantum implementation aspects, including qubit number scalings and state preparation and block-encoding costs. Section \ref{sec:Lorenz} validates the method on the Lorenz system by showing that the RBA reproduces the corresponding fully discrete Euler trajectory. Section \ref{sec:Burgers} then extends the construction to the one-dimensional viscous Burgers equation. Finally, Section \ref{sec:conclusion} summarizes the algorithm, highlighting its advantages and main limitations.

\section{reduced basis Algorithm}
\label{sec:RBA}

We first present the RBA algorithm for a system of nonlinear polynomial ODEs in Section \ref{subsec:RBA-ODEs}, and then explain how the same construction extends to nonlinear polynomial PDEs in Section \ref{subsec:RBA-PDEs}.

\subsection{RBA applied to nonlinear ODEs}
\label{subsec:RBA-ODEs}

In this section, we describe the RBA algorithm for a system of $n$ first-order\footnote{Higher-order ODEs can be rewritten as equivalent first-order systems by introducing auxiliary variables, so the same construction applies in that setting.} polynomial ODEs. We consider the initial-value problem
\begin{equation}
\frac{d\mathbf{z}(t)}{dt}
=
\mathbf{F}\bigl(\mathbf{z}(t)\bigr),
\qquad
\mathbf{z}(0)=\mathbf{b},
\label{eq:ODE-system}
\end{equation}
where $\mathbf{z}(t)=(z_1(t),\ldots,z_n(t))^T\in\mathbb{R}^n$ defines the state vector and $\mathbf{F}:\mathbb{R}^n\to\mathbb{R}^n$ is a vector-valued polynomial map with
$\mathbf{F}(\mathbf{z}(t))=(f_1(\mathbf{z}(t)),\ldots,f_n(\mathbf{z}(t)))^T$,
where each polynomial $f_\alpha$ has degree at most $p\ge 1$.

Applying the explicit Euler method\footnote{The RBA construction is not restricted to explicit Euler. Other time discretization schemes may be used, provided the resulting fully discrete update can be represented in the chosen lifted basis. We use explicit Euler here because it leads to the simplest polynomial update map and makes the construction of the RBA algorithm more transparent.} with timestep $h>0$ to Eq. \eqref{eq:ODE-system} gives the discrete time system:
\begin{equation}
\label{eq:discrete-ODE-map}
\mathbf{z}^{(k+1)}
=
\Phi_h\bigl(\mathbf{z}^{(k)}\bigr)
=
\mathbf{z}^{(k)} + h\,\mathbf{F}\bigl(\mathbf{z}^{(k)}\bigr),
\qquad
\mathbf{z}^{(0)}=\mathbf{b}.
\end{equation}
After $m$ timesteps, the discrete state generated from the initial condition is
\begin{equation}
\label{eq:composed-mapODE}
\mathbf{z}^{(m)}
=
\Phi_h^{\circ m}(\mathbf{b}),
\end{equation}
where $\Phi_h^{\circ m}$ denotes the $m$-fold composition of the Euler update map $\Phi_h$ with itself.

Since $p\ge 1$ and every component $f_\alpha$ has degree at most $p$, the Euler update map $\Phi_h$ also has polynomial components of degree at most $p$. It follows that every component of the $m$-timestep composition $\Phi_h^{\circ m}$ is a polynomial of degree at most $p^m$. To represent this nonlinear $m$-timestep composition in linear form, we lift the initial condition to a monomial basis. First, consider the full set of monomials\footnote{The constant monomial is included by allowing the multi-index $\nu=\mathbf{0}$, for which $\mathbf{z}^{\nu}=z_1^0\cdots z_n^0=1$.} in the state variables $z_1,\ldots,z_n$ of total degree at most $p^m$:
\begin{equation}
\label{eq:full-monomial-basis-ODE}
\mathcal{M}_{p^m}^{\mathrm{full}}
=
\left\{
\mathbf{z}^{\nu}
:\nu\in\mathbb{N}_0^n,\ |\nu|\le p^m
\right\},
\qquad
\mathbf{z}^{\nu}:=\prod_{i=1}^n z_i^{\nu_i}.
\end{equation}

Next, let $\psi_j^{\mathrm{red}}$ denote the individual monomials in the reduced basis
\begin{equation*}
\mathcal{M}_{p^m}^{\mathrm{red}}
=
\{\psi_1^{\mathrm{red}},\ldots,\psi_{d_{p^m}^{\mathrm{red}}}^{\mathrm{red}}\}
\subseteq
\mathcal{M}_{p^m}^{\mathrm{full}}.
\end{equation*}
This basis consists only of those monomials that appear with nonzero coefficient in at least one component of $\Phi_h^{\circ m}$. The dimension of the reduced basis satisfies the upper bound
\begin{equation}
\label{eq:ODE-reduced-dim}
d_{p^m}^{\mathrm{red}}
\le
\bigl\lvert\mathcal{M}_{p^m}^{\mathrm{full}}\bigr\rvert
=
\binom{n+p^m}{p^m}.
\end{equation}
The corresponding lifted monomial vector is
\begin{equation}
\label{eq:ODE-lifted-vector}
w_m^{\mathrm{red}}(\mathbf{z})
=
\bigl(
\psi_1^{\mathrm{red}}(\mathbf{z}),
\ldots,
\psi_{d_{p^m}^{\mathrm{red}}}^{\mathrm{red}}(\mathbf{z})
\bigr)^T.
\end{equation}

Since each component of $\Phi_h^{\circ m}(\mathbf{z})$ is a polynomial in the state variables, and since $\mathcal{M}_{p^m}^{\mathrm{red}}$ contains all monomials appearing in these polynomials, each component admits a unique expansion in the reduced  basis:
\begin{equation*}
\label{eq:ODE-polynomialbasis}
\bigl(\Phi_h^{\circ m}(\mathbf{z})\bigr)_\alpha
=
\sum_{j=1}^{d_{p^m}^{\mathrm{red}}} c_{\alpha j}^{(m)}\,\psi_j^{\mathrm{red}}(\mathbf{z}),
\qquad \alpha=1,\ldots,n.
\end{equation*}
We encode these coefficients in a padded linear operator $A_m^{\mathrm{red}}\in\mathbb{R}^{d_{p^m}^{\mathrm{red}}\times d_{p^m}^{\mathrm{red}}}$,
which we refer to as the RBA operator. Its first $n$ rows reconstruct the
physical variables after $m$ timesteps, while the remaining rows are set to zero:
\begin{equation}
\label{eq:Am-ODE-entries}
(A_m^{\mathrm{red}})_{\alpha j}
=
\begin{cases}
c_{\alpha j}^{(m)}, & 1\le \alpha \le n,\\
0, & n < \alpha \le d_{p^m}^{\mathrm{red}},
\end{cases}
\qquad j=1,\ldots,d_{p^m}^{\mathrm{red}} .
\end{equation}
Applying this operator to the lifted monomial vector defined in Eq.~\eqref{eq:ODE-lifted-vector} gives
\begin{equation}
\label{eq:Amwm_ODE-generic}
A_m^{\mathrm{red}}\,w_m^{\mathrm{red}}(\mathbf{z})
=
\begin{pmatrix}
\Phi_h^{\circ m}(\mathbf{z})\\
\mathbf{0}_{d_{p^m}^{\mathrm{red}}-n}
\end{pmatrix}.
\end{equation}
The evolved state after $m$ timesteps can therefore be obtained by applying the linear operator to the lifted monomial vector evaluated at the initial condition:
\begin{equation}
\label{eq:Amwm-ODE-initial}
A_m^{\mathrm{red}}\,w_m^{\mathrm{red}}(\mathbf{b})
=
\begin{pmatrix}
\mathbf{z}^{(m)}\\
\mathbf{0}_{d_{p^m}^{\mathrm{red}}-n}
\end{pmatrix}.
\end{equation}

This shows that the RBA construction implements the nonlinear $m$-timestep Euler map $\Phi_h^{\circ m}$ through a linear representation. Although the map is nonlinear in the original state variables $\mathbf{z}$, its action is recovered by applying the linear operator $A_m^{\mathrm{red}}$ to the lifted monomial vector $w_m^{\mathrm{red}}(\mathbf{b})$. Thus, the RBA representation is exact at the level of the discrete Euler dynamics,
and introduces no additional approximation error beyond the Euler time discretization error.

\subsection{RBA applied to nonlinear PDEs}
\label{subsec:RBA-PDEs}

The RBA algorithm extends naturally to polynomial nonlinear PDEs after spatial discretization. Consider a system of $n$ coupled evolution equations on a spatial domain 
$\Omega\subseteq\mathbb{R}^D$,
\begin{equation}
\label{eq:PDE-system}
\begin{cases}
\partial_t \mathbf{z}(x,t)
=
\mathcal{N}\bigl(\mathbf{z}\bigr)(x,t),
& (x,t)\in \Omega\times(0,T],\\[1mm]
\mathbf{z}(x,0)=\mathbf{z}_0(x),
& x\in\Omega,
\end{cases}
\end{equation}
together with appropriate boundary conditions on $\partial\Omega\times(0,T]$.
Here
$\mathbf{z}(x,t)=(z_1(x,t),\ldots,z_n(x,t))^T\in\mathbb{R}^n$
denotes the vector of coupled fields, and $\mathcal{N}$ is a possibly nonlinear
differential operator.

We discretize $\Omega$ using $N$ degrees of freedom in each of the $D$ spatial
directions. Let
$\boldsymbol{\ell}=(\ell_1,\ldots,\ell_D)\in\{1,\ldots,N\}^D$
denote a spatial multi-index, and let $z_{j,\boldsymbol{\ell}}(t)$ denote the
discrete value of the $j$-th field at a grid point indexed by
$\boldsymbol{\ell}$. Collecting all discrete unknowns into a single
semi-discrete state vector gives
\begin{equation*}
\mathbf{z}(t)
=
\bigl(z_{j,\boldsymbol{\ell}}(t)\bigr)^T_
{\substack{1\le j\le n\\ \boldsymbol{\ell}\in\{1,\ldots,N\}^D}}
\in
\mathbb{R}^{nN^D}.
\end{equation*}

Applying the method of lines to Eq. \eqref{eq:PDE-system}, one obtains the finite-dimensional initial-value problem
\begin{equation}
\label{eq:PDEtoODE}
\frac{d\mathbf{z}(t)}{dt}
=
\mathbf{F}\bigl(\mathbf{z}(t)\bigr),
\qquad
\mathbf{z}(0)=\mathbf{b},
\end{equation}
where $\mathbf{b}\in\mathbb{R}^{nN^D}$ is the discretized initial condition. The boundary conditions are incorporated into the definition of the semi-discrete vector field $\mathbf{F}:\mathbb{R}^{nN^D}\to\mathbb{R}^{nN^D}$. We assume that the spatial discretization of $\mathcal{N}$ produces a polynomial vector field in the state variables. The resulting semi-discrete vector field $\mathbf{F}$ is a vector-valued polynomial map, $\mathbf{F}(\mathbf{z}) = \bigl(f_1(\mathbf{z}),\ldots,f_{nN^D}(\mathbf{z})\bigr)^T$, where each polynomial $f_\alpha$ has degree at most $p\ge 1$.

At this point, we apply to Eq. \eqref{eq:PDEtoODE} the same explicit Euler discretization used in the ODE construction. This gives the same discrete time
polynomial update map defined in Eq. \eqref{eq:discrete-ODE-map}, and the
$m$-timestep discrete state generated by polynomial composition as in Eq. \eqref{eq:composed-mapODE}. The only change is that the state dimension is now $nN^D$ rather than $n$. Consequently, by the same degree-counting argument as in the ODE case, every component of the composed map $\Phi_h^{\circ m}$ is a polynomial of degree at most $p^m$.

We now apply the same lifting construction as in Section \ref{subsec:RBA-ODEs} by first considering the full global monomial basis
\begin{equation}
\label{eq:PDE-global-full basis}
\mathcal{M}_{p^m}^{\mathrm{glob,full}}
=
\left\{
\mathbf{z}^{\nu}
:\nu\in\mathbb{N}_0^{nN^D},\ |\nu|\le p^m
\right\},
\qquad
\mathbf{z}^{\nu}:=\prod_{i=1}^{nN^D} z_i^{\nu_i}.
\end{equation}
We call this basis global because it is defined over the entire semi-discrete state vector. Its dimension is
\begin{equation*}
\label{eq:PDE-global-full-dim}
\bigl|\mathcal{M}_{p^m}^{\mathrm{glob,full}}\bigr|
=
\binom{nN^D+p^m}{p^m}.
\end{equation*}

Working with the global full monomial basis is usually far too pessimistic for PDE discretizations, since it ignores spatial locality. For many finite-difference, finite-volume, and finite-element discretizations, the value
at a grid point $\boldsymbol{\ell}$ after $m$ timesteps depends only on the field values in an effective stencil
$S_m(\boldsymbol{\ell})\subseteq \{1,\ldots,N\}^D$.
Let 
\begin{equation*} 
s_m(\boldsymbol{\ell}) = \bigl\lvert S_m(\boldsymbol{\ell})\bigr\rvert
\end{equation*}
denote the size of this stencil. Since each grid point carries $n$ field
components, the number of scalar variables in the corresponding local stencil is
\begin{equation*}
\label{eq:nr-scalar-variable}
v_m(\boldsymbol{\ell})
=
n\,s_m(\boldsymbol{\ell}).
\end{equation*}

We collect these variables into a local stencil vector $\mathbf{x}^{(\boldsymbol{\ell})}\in\mathbb{R}^{v_m(\boldsymbol{\ell})}$. Then, for each field component $\alpha=1,\ldots,n$, the component value $\bigl(\Phi_h^{\circ m}(\mathbf{z})\bigr)_{\alpha,\boldsymbol{\ell}}$
depends only on $\mathbf{x}^{(\boldsymbol{\ell})}$. Hence, there exists a polynomial
\begin{equation*}
\label{eq:local-poly}
g_{\alpha,\boldsymbol{\ell}}^{(m)}
:
\mathbb{R}^{v_m(\boldsymbol{\ell})}
\to
\mathbb{R},
\qquad
\deg g_{\alpha,\boldsymbol{\ell}}^{(m)}\le p^m,
\end{equation*}
such that
\begin{equation*}
\label{eq:local-m-step}
\bigl(\Phi_h^{\circ m}(\mathbf{z})\bigr)_{\alpha,\boldsymbol{\ell}}
=
g_{\alpha,\boldsymbol{\ell}}^{(m)}
\bigl(\mathbf{x}^{(\boldsymbol{\ell})}\bigr).
\end{equation*}
Thus, the exact $m$-timestep update at a grid point can be represented using only the variables in its effective stencil, rather than all $nN^D$ variables in the global semi-discrete state vector. 

We therefore introduce the full local monomial basis
\begin{equation*}
\label{eq:local-full-monomial-basis}
\mathcal{M}_{p^m}^{\mathrm{loc,full}}(\boldsymbol{\ell})
=
\left\{
\bigl(\mathbf{x}^{(\boldsymbol{\ell})}\bigr)^\nu
:
\nu\in\mathbb{N}_0^{v_m(\boldsymbol{\ell})},
\ |\nu|\le p^m
\right\},
\end{equation*}
where
\begin{equation*}
\bigl(\mathbf{x}^{(\boldsymbol{\ell})}\bigr)^\nu
:=
\prod_{i=1}^{v_m(\boldsymbol{\ell})}
\left(x_i^{(\boldsymbol{\ell})}\right)^{\nu_i}.
\end{equation*}
Let $\psi_{j,\boldsymbol{\ell}}^{\mathrm{loc},red}$ denote the individual monomials in the reduced local monomial basis
\begin{equation}
\label{eq:monomials-mlocred}
\mathcal{M}_{p^m}^{\mathrm{loc,red}}(\boldsymbol{\ell}) = \left\{
    \psi_{1,\boldsymbol{\ell}}^{\mathrm{loc},red},
    \ldots,
    \psi_{d_{p^m}^{\mathrm{loc,red}}(\boldsymbol{\ell}),\boldsymbol{\ell}}^{\mathrm{loc,red}}
    \right\}
\subseteq
\mathcal{M}_{p^m}^{\mathrm{loc,full}}(\boldsymbol{\ell}),
\end{equation}
which contains only the monomials that appear with nonzero coefficient in at least one of the polynomials
$g_{\alpha,\boldsymbol{\ell}}^{(m)}$, for $\alpha=1,\ldots,n$. The size of this basis satisfies the following upper bound 
\begin{equation}
\label{eq:dloc-PDE}
d_{p^m}^{\mathrm{loc,red}}(\boldsymbol{\ell})
\le
\left|
\mathcal{M}_{p^m}^{\mathrm{loc,full}}(\boldsymbol{\ell})
\right|
=
\binom{v_m(\boldsymbol{\ell})+p^m}{p^m}.
\end{equation}

The corresponding local lifted monomial vector is
\begin{equation}
\label{eq:pde-local-lifted-state}
w_{m,\boldsymbol{\ell}}^{\mathrm{loc,red}}
\bigl(\mathbf{x}^{(\boldsymbol{\ell})}\bigr)
=
\left(
\psi_{1,\boldsymbol{\ell}}^{\mathrm{loc,red}}
\bigl(\mathbf{x}^{(\boldsymbol{\ell})}\bigr),
\ldots,
\psi_{d_{p^m}^{\mathrm{loc,red}}(\boldsymbol{\ell}),\boldsymbol{\ell}}^{\mathrm{loc,red}}
\bigl(\mathbf{x}^{(\boldsymbol{\ell})}\bigr)
\right)^T .
\end{equation}

Since the reduced local basis contains all monomials appearing in the local
update polynomials, each component admits a unique expansion
\begin{equation*}
\label{eq:local-expansion}
g_{\alpha,\boldsymbol{\ell}}^{(m)}
\bigl(\mathbf{x}^{(\boldsymbol{\ell})}\bigr)
=
\sum_{j=1}^{d_{p^m}^{\mathrm{loc,red}}(\boldsymbol{\ell})}
c_{\alpha j}^{(m,\boldsymbol{\ell})}
\,
\psi_{j,\boldsymbol{\ell}}^{\mathrm{loc,red}}
\bigl(\mathbf{x}^{(\boldsymbol{\ell})}\bigr),
\qquad
\alpha=1,\ldots,n.
\end{equation*}
We encode these coefficients in a padded local linear operator $A_{m,\boldsymbol{\ell}}^{\mathrm{loc,red}}
\in
\mathbb{R}^{
d_{p^m}^{\mathrm{loc,red}}(\boldsymbol{\ell})
\times
d_{p^m}^{\mathrm{loc,red}}(\boldsymbol{\ell})
}$, whose first $n$ rows reconstruct the $n$ field components at grid point
$\boldsymbol{\ell}$ after $m$ timesteps, while the remaining rows are set to
zero. That is, for
$\alpha,j=1,\ldots,d_{p^m}^{\mathrm{loc,red}}(\boldsymbol{\ell})$,
\begin{equation}
\label{eq:local-Am-entries}
\bigl(A_{m,\boldsymbol{\ell}}^{\mathrm{loc,red}}\bigr)_{\alpha j}
=
\begin{cases}
c_{\alpha j}^{(m,\boldsymbol{\ell})},
& 1\le \alpha \le n,\\
0,
& n<\alpha.
\end{cases}
\end{equation}
Applying this operator to the lifted monomial vector defined in Eq. \eqref{eq:pde-local-lifted-state} gives
\begin{equation*}
\label{eq:local-Am-action}
A_{m,\boldsymbol{\ell}}^{\mathrm{loc,red}}
\,
w_{m,\boldsymbol{\ell}}^{\mathrm{loc,red}}
\bigl(\mathbf{x}^{(\boldsymbol{\ell})}\bigr)
=
\begin{pmatrix}
g_{1,\boldsymbol{\ell}}^{(m)}
\bigl(\mathbf{x}^{(\boldsymbol{\ell})}\bigr)\\
\vdots\\
g_{n,\boldsymbol{\ell}}^{(m)}
\bigl(\mathbf{x}^{(\boldsymbol{\ell})}\bigr)\\
\mathbf{0}_{d_{p^m}^{\mathrm{loc,red}}(\boldsymbol{\ell})-n}
\end{pmatrix}.
\end{equation*}
For a prescribed discretized initial condition $\mathbf{z}^{(0)}=\mathbf{b}$, we denote by $\mathbf{x}^{(\boldsymbol{\ell})}(\mathbf{b})$ the restriction of $\mathbf{b}$ to the scalar variables in the effective stencil $S_m(\boldsymbol{\ell})$. Evaluating the local lifted vector at this restricted initial condition gives
\begin{equation*}
\label{eq:local-Am-action-initial}
A_{m,\boldsymbol{\ell}}^{\mathrm{loc,red}}
\,
w_{m,\boldsymbol{\ell}}^{\mathrm{loc,red}}
\bigl(\mathbf{x}^{(\boldsymbol{\ell})}(\mathbf{b})\bigr)
=
\begin{pmatrix}
g_{1,\boldsymbol{\ell}}^{(m)}
\bigl(\mathbf{x}^{(\boldsymbol{\ell})}(\mathbf{b})\bigr)\\
\vdots\\
g_{n,\boldsymbol{\ell}}^{(m)}
\bigl(\mathbf{x}^{(\boldsymbol{\ell})}(\mathbf{b})\bigr)\\
\mathbf{0}_{d_{p^m}^{\mathrm{loc,red}}(\boldsymbol{\ell})-n}
\end{pmatrix}
=
\begin{pmatrix}
z_{1,\boldsymbol{\ell}}^{(m)}\\
\vdots\\
z_{n,\boldsymbol{\ell}}^{(m)}\\
\mathbf{0}_{d_{p^m}^{\mathrm{loc,red}}(\boldsymbol{\ell})-n}
\end{pmatrix}.
\end{equation*}
Thus, the nonlinear $m$-timestep update at a grid point is represented exactly by a linear operator acting on a local lifted monomial vector evaluated on the corresponding local stencil of the initial condition.

To combine the local constructions over the full grid, we introduce a
position register indexed by the spatial multi-index $\boldsymbol{\ell}$,
with basis states $\ket{\boldsymbol{\ell}}$ for
$\boldsymbol{\ell}\in\{1,\ldots,N\}^D$. Since the reduced local lifted
dimension $d_{p^m}^{\mathrm{loc,red}}(\boldsymbol{\ell})$ may depend on
the grid point, for example because interior and boundary grid points
can have different effective stencils, we embed all local lifted spaces
into a common local register of dimension
\begin{equation}
\label{eq:PDE-common-local-dimension}
d_{p^m}^{\mathrm{loc,max}}
=
\max_{\boldsymbol{\ell}\in\{1,\ldots,N\}^D}
d_{p^m}^{\mathrm{loc,red}}(\boldsymbol{\ell}).
\end{equation}
The local lifted vector
$w_{m,\boldsymbol{\ell}}^{\mathrm{loc,red}}
(\mathbf{x}^{(\boldsymbol{\ell})})$ is then zero-padded to define
\begin{equation}
\label{eq:paddedlorredvector}
\widetilde{w}_{m,\boldsymbol{\ell}}^{\mathrm{loc,red}}
(\mathbf{x}^{(\boldsymbol{\ell})})
\in
\mathbb{R}^{d_{p^m}^{\mathrm{loc,max}}}.
\end{equation}
Equivalently, the local operator
$A_{m,\boldsymbol{\ell}}^{\mathrm{loc,red}}$ is extended by zero to an
operator
\begin{equation}
\label{eq:AMred-padded}
    \widetilde{A}_{m,\boldsymbol{\ell}}^{\mathrm{loc,red}}
    :
    \mathbb{R}^{d_{p^m}^{\mathrm{loc,max}}}
    \to
    \mathbb{R}^{d_{p^m}^{\mathrm{loc,max}}},
\end{equation}
acting on the common local lifted space. With this convention, the
position-controlled local RBA operator is
\begin{equation}
\label{eq:position-controlled-Am}
\mathcal{A}_m^{\mathrm{red}}
=
\sum_{\boldsymbol{\ell}\in\{1,\ldots,N\}^D}
\ket{\boldsymbol{\ell}}\!\bra{\boldsymbol{\ell}}
\otimes
\widetilde{A}_{m,\boldsymbol{\ell}}^{\mathrm{loc,red}} .
\end{equation}

\section{Quantum Algorithm Construction}
\label{sec:QAC}

As the RBA algorithm is a general algorithm that can be applied to a wide array of polynomial nonlinear differential equations, its implementation on a quantum computer depends on the specific problem it is used to solve. Therefore, it is difficult to estimate the cost of running such an algorithm on quantum hardware, since this cost is highly problem-specific. Nevertheless, in this section we discuss several aspects that need to be considered when estimating the algorithmic cost for a particular problem. These include the number of qubits required to encode and process the lifted quantum state, the cost of preparing this state, the classical preprocessing cost required to compute the reduced basis and the RBA operator, and the question of how to efficiently block-encode this operator given its structure.

\subsection{Quantum state encoding}

\paragraph{Qubit number scaling.}For a general system of $n$ nonlinear ODEs, the lifted quantum state associated with the $m$-timestep reduced representation is
\begin{equation}
\label{eq:ODE-amplitude-encoding}
\ket{w_m(\mathbf{z})}
=
\frac{1}{\|w_m^{\mathrm{red}}(\mathbf{z})\|_2}
\sum_{j=1}^{d_{p^m}^{\mathrm{red}}}
\bigl(w_m^{\mathrm{red}}(\mathbf{z})\bigr)_j
\ket{j}.
\end{equation}
Here, $w_m^{\mathrm{red}}(\mathbf{z})$ is the lifted monomial vector constructed from the reduced ordered basis $\mathcal{M}_{p^m}^{\mathrm{red}}$, as defined in Eq. \eqref{eq:ODE-lifted-vector}, and the index $j$ labels the $j$-th monomial in this basis. Ignoring ancilla qubits required for block-encoding, the number of qubits needed to encode the lifted ODE register is
\begin{equation*}
\label{eq:ODE-qubit-count-1}
q_m^{\mathrm{ODE}}
=
\left\lceil \log_2 d_{p^m}^{\mathrm{red}} \right\rceil .
\end{equation*}

Since the dimension of the reduced basis $d_{p^m}^{\mathrm{red}}$ satisfies Eq. \eqref{eq:ODE-reduced-dim}, we obtain the following upper bound for the number of qubits needed to encode the full basis
\begin{equation}
\label{eq:ODE-qubit-count}
q_m^{\mathrm{ODE}}
\le
\left\lceil
\log_2
\binom{n+p^m}{p^m}
\right\rceil .
\end{equation}
For $p>1$, this bound scales as 
\begin{equation*}
\left\lceil
\log_2
\binom{n+p^m}{p^m}
\right\rceil  = O(nm\log p),
\end{equation*}
since
\begin{equation*}
\binom{n+p^m}{p^m}=\binom{n+p^m}{n}=O((p^m)^n).
\end{equation*}
Hence, the full lifted basis dimension grows exponentially with $m$, but the number of qubits needed to index the lifted basis grows at most linearly in $m$, with
\begin{equation*}
    \label{eq:ODE-qubits}
    q_m^{\mathrm{ODE}}=O(nm\log p).
\end{equation*}
This is a worst-case estimate since the actual reduced basis qubit count can be smaller because $\mathcal{M}_{p^m}^{\mathrm{red}}$ contains only the monomials appearing with nonzero coefficient in the composed $m$-timestep map.

For spatially discretized PDEs, the lifted quantum state contains both a
position register $\ket{\boldsymbol{\ell}}$ and a local lifted basis
register $\ket{j}$. Using the common local dimension
$d_{p^m}^{\mathrm{loc,max}}$ introduced in
Eq. \eqref{eq:PDE-common-local-dimension}, the global lifted vector can
be written in block form as
\begin{equation*}
\label{eq:PDE-lifted-block-vector}
W_m(\mathbf{z})
=
\sum_{\boldsymbol{\ell}\in\{1,\ldots,N\}^D}
\ket{\boldsymbol{\ell}}
\otimes
\widetilde{w}_{m,\boldsymbol{\ell}}^{\mathrm{loc,red}}
\bigl(\mathbf{x}^{(\boldsymbol{\ell})}\bigr).
\end{equation*}
Its amplitude-encoded quantum state is
\begin{equation}
\label{eq:PDE-amplitude-encoding}
\ket{W_m(\mathbf{z})}
=
\frac{1}{\|W_m(\mathbf{z})\|_2}
\sum_{\boldsymbol{\ell}\in\{1,\ldots,N\}^D}
\sum_{j=1}^{d_{p^m}^{\mathrm{loc,max}}}
\bigl(
\widetilde{w}_{m,\boldsymbol{\ell}}^{\mathrm{loc,red}}
(\mathbf{x}^{(\boldsymbol{\ell})})
\bigr)_j
\ket{\boldsymbol{\ell}}\ket{j}.
\end{equation}
Here, $\widetilde{w}_{m,\boldsymbol{\ell}}^{\mathrm{loc,red}}$ denotes the zero-padded version of the local lifted monomial vector
$w_{m,\boldsymbol{\ell}}^{\mathrm{loc,red}}$ defined in Eq. \eqref{eq:pde-local-lifted-state}.

Ignoring ancilla qubits required for block-encoding, the number of qubits
needed for the lifted PDE representation is
\begin{equation*}
\label{eq:PDE-qubit-count}
q_m^{\mathrm{PDE}}
=
\left\lceil \log_2(N^D)\right\rceil
+
\left\lceil \log_2 d_{p^m}^{\mathrm{loc,max}}\right\rceil .
\end{equation*}
The first term indexes the spatial grid point
$\boldsymbol{\ell}\in\{1,\ldots,N\}^D$, while the second term indexes the
common local lifted basis register. Since each reduced local basis is a subset of the corresponding full local
monomial basis, we have
\begin{equation*}
d_{p^m}^{\mathrm{loc,max}}
\le
\binom{v_m^{\max}+p^m}{p^m},
\qquad
v_m^{\max}
=
\max_{\boldsymbol{\ell}\in\{1,\ldots,N\}^D}
v_m(\boldsymbol{\ell}),
\end{equation*}
where $v_m^{\max}$ is the maximal number of scalar variables appearing in
an $m$-timestep effective stencil. Therefore,
\begin{equation}
\label{eq:PDE-qubit-count-bound}
q_m^{\mathrm{PDE}}
\le
\left\lceil \log_2(N^D)\right\rceil
+
\left\lceil
\log_2
\binom{v_m^{\max}+p^m}{p^m}
\right\rceil .
\end{equation}

For a local discretization in $D$ spatial dimensions, the position register contributes
\begin{equation*}
\label{eq:grid_qubits_PDE}
    \left\lceil \log_2(N^D)\right\rceil
    =
    O(D\log N)
\end{equation*}
qubits. The remaining contribution comes from the local lifted register. 
Assume that the spatial discretization is explicit and local, with a fixed
one-timestep stencil radius independent of $N$. Then, after $m$ composed
timesteps, the effective stencil can expand by at most $O(m)$ grid points in
each spatial direction, provided that this propagated stencil remains small
compared with the full spatial domain. Hence, for $m \ll N$, up to constants
depending on the one-timestep stencil radius, the number of grid points in an
$m$-timestep effective stencil satisfies
\begin{equation*}
        s_m^{\max}=O(m^D).
\end{equation*}

Since each grid point carries $n$ scalar field variables, this gives
\begin{equation*}
        v_m^{\max}=O(nm^D).
\end{equation*}
Therefore,
\begin{equation*}
\log_2
\binom{v_m^{\max}+p^m}{p^m}
=
\log_2
\binom{v_m^{\max}+p^m}{v_m^{\max}} .
\end{equation*}
For $p>1$, this term can be bounded as
\begin{equation*}
\log_2
\binom{v_m^{\max}+p^m}{v_m^{\max}}
=
O\!\left(v_m^{\max}m\log p\right)
=
O\!\left(nm^{D+1}\log p\right).
\end{equation*}
Consequently, under the locality assumptions above and as long as the
$m$-timestep effective stencil remains small compared with the full spatial
domain, the PDE qubit count satisfies
\begin{equation*}
      \label{eq:PDE-qubits}
    q_m^{\mathrm{PDE}}
    =
    O(D\log N+n m^{D+1}\log p).
\end{equation*}

Thus, in this regime, the position register contributes only logarithmically in the number of grid points per spatial direction, while the nonlinear lifted register overhead is controlled by the size of the propagated local stencil and by the polynomial degree $p^m$. If $m$ becomes large enough that the effective stencil is comparable to the full spatial domain, then the locality assumption no longer gives an $N$-independent bound for the local register. In that case, the local stencil size must be replaced by its saturated value, which may be as large as $O(N^D)$, and the bound above no longer applies.

This polynomial-in-time qubit scaling improves on the copy-based approach of Leyton and Osborne\cite{leyton2008quantum}, where the number of qubits required to run the algorithm grows exponentially in time. It is more comparable to the polynomial time dependence obtained in the mean-field approach\cite{lloyd2020quantum}, which avoids the exponential copy overhead by using many weakly interacting copies. However, the mechanism is different. The mean-field method approximates the nonlinear dynamics with an error controlled by the number of copies and by stability properties of the evolution. In contrast, the RBA construction represents the $m$-timestep discrete nonlinear dynamics exactly, with no additional closure, truncation, or mean-field error introduced.

\paragraph{State preparation cost.}
The quantum implementation requires preparing the amplitude-encoded lifted state, namely $\lvert{w_m(\mathbf{z})}\rangle$ in the ODE case and $\lvert{W_m(\mathbf{z})}\rangle$ in the PDE case, as defined in Eqs. \eqref{eq:ODE-amplitude-encoding} and
\eqref{eq:PDE-amplitude-encoding}, respectively. A conservative estimate of the state preparation cost is obtained by treating the lifted monomial vector as an arbitrary amplitude vector. For the ODE register, which contains $q_m^{\mathrm{ODE}}$ qubits as defined in Eq. \eqref{eq:ODE-qubit-count}, the generic state preparation result of Li et al. \cite{li2026reducing} gives an upper
bound of
\begin{equation*}
\frac{11}{12}2^{q_m^{\mathrm{ODE}}}
\end{equation*}
CNOT gates for preparing this state. Thus, in the worst case, the state preparation cost is linear in the dimension
of the reduced lifted space, since
\begin{equation*}
\label{eq:ODE-stateprep}
2^{q_m^{\mathrm{ODE}}}
=
O\!\left(d_{p^m}^{\mathrm{red}}\right).
\end{equation*}

In the PDE setting, the same arbitrary state preparation bound applied to the full lifted PDE register gives
\begin{equation*}
\frac{11}{12}2^{q_m^{\mathrm{PDE}}}
\end{equation*}
CNOT gates, where $q_m^{\mathrm{PDE}}$ is defined in
Eq. \eqref{eq:PDE-qubit-count-bound}. Equivalently, this worst-case cost is linear in the product of the number of grid points and the padded local lifted dimension:
\begin{equation*}
\label{eq:PDE-stateprep}
2^{q_m^{\mathrm{PDE}}}
=
O\!\left(N^D d_{p^m}^{\mathrm{loc,max}}\right).
\end{equation*}

However, one must keep in mind that this state preparation estimate is a generic worst-case estimate. The lifted vectors constructed in the RBA algorithm are not arbitrary amplitude vectors. The lifted states inherit algebraic structure from the monomials, the sparsity of the reduced basis, and, in the PDE case, the repeated local stencil structure. Nevertheless, the existence of such structure does not by itself imply that the corresponding quantum state can be prepared efficiently. To improve on the generic cost, one would need to show that this structure can be exploited by known simple quantum operations or by a problem-specific state preparation procedure. This question is problem dependent and lies outside the scope of the present work. For this reason, the bounds presented in this section should be understood as generic worst-case upper bounds for RBA state preparation in the absence of an explicit efficient state preparation circuit.

\subsection{RBA Operator Construction and Block-Encoding}
\label{subsec:pre-processing}

\paragraph{Classical preprocessing stage.} To run the RBA algorithm on a quantum computer, we must first perform a
classical preprocessing step in which we construct the RBA operator required for the $m$-timestep update. For the ODE case, this is the operator $A_m^{\mathrm{red}}$ defined in Eq. \eqref{eq:Am-ODE-entries}. For the PDE case, this is the operator $\mathcal{A}_m$ defined in Eq. \eqref{eq:position-controlled-Am}, or more specifically the reduced local operators $A_{m,\boldsymbol{\ell}}^{\mathrm{loc,red}}$ defined in Eq. \eqref{eq:local-Am-entries}. In a second step, these operators must be block-encoded such that they can admit a unitary implementation.

The classical preprocessing consists of two main tasks. First, we must identify the reduced monomial basis appearing in the composed $m$-timestep polynomial map $\Phi_h^{\circ m}$. Second, we must construct the corresponding RBA operator. Naturally, the cost of this step is problem specific. It depends on the sparse polynomial supports that are generated during the construction of $\Phi_h^{\circ m}$, and not directly
on the size of the full monomial basis, unless the reduced basis coincides with the full basis. In many PDE discretizations, the full basis case is overly pessimistic. The
local stencil, the sparsity of the nonlinear terms, and the structure of the discretization can substantially reduce the number of monomials that actually appear in the composed map.

More precisely, the reduced basis does not need to be constructed by looping over the full monomial set $\mathcal{M}_{p^m}^{\mathrm{full}}$. Instead, we represent each
polynomial component sparsely, using exponent vectors as keys and the corresponding coefficients as values. Starting from the one-timestep map $\Phi_h$, we form the composed maps
\begin{equation*}
P^{(k)}=\Phi_h^{\circ k}, \qquad k=1,\ldots,m,
\end{equation*}
by sparse substitution. At each composition step, each monomial $z^\nu$ appearing in a component of $\Phi_h$ is replaced by
\begin{equation*}
z^\nu \mapsto \prod_i \left(P_i^{(k)}(z)\right)^{\nu_i},
\end{equation*}
where the products are evaluated using sparse polynomial multiplication and terms with the same monomial are combined. After $P^{(m)}$ has been constructed, the reduced basis is obtained as
\begin{equation*}
\mathcal{M}_{p^m}^{\mathrm{red}}
=
\bigcup_{\alpha=1}^{n} \operatorname{supp} P_\alpha^{(m)}.
\end{equation*}
Thus, the construction only generates monomials that arise through the sparse composition process. It does not test all monomials in $\mathcal{M}_{p^m}^{\mathrm{full}}$ and reject those with zero coefficient.

This distinction is important for the preprocessing cost. Let
\begin{equation*}
\Sigma_k=\sum_{\alpha=1}^n |\operatorname{supp} P_\alpha^{(k)}|
\end{equation*}
denote the total sparse support size of the $k$-timestep map, and let $G_k$ denote the number of monomial products generated while forming $P^{(k+1)}=\Phi_h\circ P^{(k)}$. With naive sparse polynomial multiplication, the preprocessing cost is controlled by
\begin{equation*}
T_{\mathrm{pre}}
=
O\!\left(\sum_{k=1}^{m-1} G_k + \Sigma_m \log \Sigma_m\right),
\end{equation*}
up to coefficient arithmetic costs. The first term accounts for the sparse polynomial compositions used to construct the $m$-timestep map, while the second term accounts for collecting and ordering the final reduced basis. Consequently, this sparse construction avoids the enumeration of the full basis. The preprocessing cost is instead controlled by the intermediate sparse products
generated during the compositions, represented by the quantities $G_k$, together with the final sparse support $\Sigma_m$. These supports may still grow exponentially with $m$, even when they are much smaller than the full
monomial basis. In such cases, the identification of the reduced basis and the construction of the corresponding RBA operator also become exponentially expensive in the timestep window $m$.

For PDEs, the local RBA operator does not need to be constructed independently at every grid point. For local discretizations on structured grids, the algebraic form of the update is typically the same at all interior points whose $m$-timestep effective stencil remains away from the boundary. Hence, a single local RBA operator can be constructed for the interior stencil and reused across all such grid points, with only the mapping from local stencil variables to global grid indices changing from point to point. Boundary-adjacent points, whose effective stencils intersect the boundary, may require distinct local RBA operators. Consequently, the preprocessing cost is governed by the number of distinct effective stencil types rather than directly by the total number of grid points.

As $m$ increases, the effective stencil radius grows. The boundary layer requiring boundary-dependent operators therefore also grows, while the region where the interior operator can be reused becomes smaller. If the effective stencil becomes comparable to the domain size, the number of distinct local RBA operators can increase, making the local construction less favorable. Thus, the RBA construction is naturally suited to short timestep windows. In practice, this means that one may apply the RBA operator over a fixed $m$-timestep window, update the physical state, and then reinitialize the lifted state using this updated physical state as the new initial condition for the next window. For fixed $m$, the reduced basis and the corresponding RBA operator depend only on the discrete update rule, the timestep size, the model parameters, and the stencil geometry, not on the current value of the state.
They can therefore be constructed once at the start of the simulation and reused at every window. The drawback is that, in a quantum implementation, each reinitialization requires preparing the lifted monomial state from the updated physical state. Unless an efficient re-lifting or state-preparation procedure is available, this may require reading out the updated state, potentially through quantum state tomography or an equivalent full-state estimation procedure.

For PDE discretizations, if one wants to understand how the local stencil develops over several timesteps, this can often be done by propagating the one-timestep stencil before carrying out the polynomial composition. The precise growth of this stencil depends on the discretization. For example, axis-aligned nearest-neighbor discretizations lead to a propagated stencil of von Neumann type, while tensor-product stencils with diagonal couplings, higher-order finite-difference stencils, or wider finite-volume reconstructions lead to different stencil geometries. In many standard cases, the size of the propagated stencil can be obtained directly from the one-timestep stencil and written in closed form, or bounded by a simple expression. This is useful for estimating the maximal local variable count entering Eq. \eqref{eq:dloc-PDE}.

However, this stencil information only determines which variables may enter the local polynomial. It does not, in general, determine which monomial products of these variables occur with nonzero coefficient. The reduced local basis must still be obtained by sparse local composition of the stencil
polynomial and by collecting the final support,
\begin{equation*}
\mathcal{M}_{p^m}^{\mathrm{loc,red}}(\ell)
=
\bigcup_{\alpha=1}^n
\operatorname{supp} g_{\alpha,\ell}^{(m)}.
\end{equation*}

\paragraph{Low-rank block-encoding.}
An important structural property of the RBA operator is that it is low rank relative to the lifted dimension. In the ODE construction, the operator $A_m^{\mathrm{red}}$ defined in Eq. \eqref{eq:Am-ODE-entries} has nonzero entries only in its first $n$ rows, since its role is to reconstruct the $n$ physical variables from the lifted monomial vector. Hence,
\begin{equation*}
    \operatorname{rank}\!\left(A_m^{\mathrm{red}}\right) \leq n .
\end{equation*}
Similarly, in the local PDE construction, each local RBA operator $A_{m,\boldsymbol{\ell}}^{\mathrm{loc,red}}$ defined in
Eq. \eqref{eq:local-Am-entries} has nonzero entries only in its first $n$ rows, corresponding to the $n$ field components at the grid point $\boldsymbol{\ell}$. Therefore,
\begin{equation*}
    \operatorname{rank}\!\left(
    A_{m,\boldsymbol{\ell}}^{\mathrm{loc,red}}
    \right)
    \leq n .
\end{equation*}

This low-rank structure can be exploited when block-encoding the RBA matrices. The low-rank block-encoding protocol of Li et al. \cite{li2026reducing} gives a single ancilla $(\|A\|_2,1)$-block-encoding of a rank-$K$ matrix. In the notation used
here, where $q$ denotes the number of system qubits and $q+1$ is the total number of qubits including the block-encoding ancilla, the leading CNOT count is
\begin{equation*}
    \left(K+\frac{11}{12}\right)2^{q+1}
    +
    O\!\left(K(q+1)^2\right).
\end{equation*}

Applying this result to the ODE setting gives the upper bound
\begin{equation*}
    N_{\mathrm{CNOT}}\!\left(A_m^{\mathrm{red}}\right)
    \le
    \left(K_m+\frac{11}{12}\right)2^{q_m^{\mathrm{ODE}}+1}
    +
    O\!\left(K_m(q_m^{\mathrm{ODE}}+1)^2\right),
\end{equation*}
where
\begin{equation*}
    K_m
    =
    \operatorname{rank}\!\left(A_m^{\mathrm{red}}\right)
    \le n .
\end{equation*}
Since $2^{q_m^{\mathrm{ODE}}}=O\!\left(d_{p^m}^{\mathrm{red}}\right)$,
the leading CNOT scaling becomes
\begin{equation*}
    N_{\mathrm{CNOT}}\!\left(A_m^{\mathrm{red}}\right)
    =
    O\!\left(K_m d_{p^m}^{\mathrm{red}}\right),
\end{equation*}
which is linear in the reduced lifted basis dimension for fixed rank
$K_m$.

For the local PDE construction, the same estimate applies to each local
RBA operator. Defining
\begin{equation*}
    K_{m,\boldsymbol{\ell}}
    =
    \operatorname{rank}\!\left(
    A_{m,\boldsymbol{\ell}}^{\mathrm{loc,red}}
    \right)
    \le n ,
\end{equation*}
we obtain
\begin{equation*}
    N_{\mathrm{CNOT}}\!\left(
    A^{\mathrm{loc,red}}_{m,\boldsymbol{\ell}}
    \right)
    \le
    \left(K_{m,\boldsymbol{\ell}}+\frac{11}{12}\right)
    2^{q_m^{\mathrm{loc}}+1}
    +
    O\!\left(
    K_{m,\boldsymbol{\ell}}(q_m^{\mathrm{loc}}+1)^2
    \right),
\end{equation*}
where $q_m^{\mathrm{loc}}$ denotes the number of qubits in the common local
lifted basis register
\begin{equation*}
    q_m^{\mathrm{loc}}
    =
    \left\lceil
    \log_2 d_{p^m}^{\mathrm{loc,max}}
    \right\rceil .
\end{equation*}
Thus, the leading local block-encoding cost satisfies
\begin{equation*}
    N_{\mathrm{CNOT}}\!\left(
    A^{\mathrm{loc,red}}_{m,\boldsymbol{\ell}}
    \right)
    =
    O\!\left(
    K_{m,\boldsymbol{\ell}}d_{p^m}^{\mathrm{loc,max}}
    \right).
\end{equation*}

The local block-encodings are then combined according to the global
position-controlled RBA operator $\mathcal{A}_m$ defined in
Eq. \eqref{eq:position-controlled-Am}. The position register
$\ket{\boldsymbol{\ell}}$ selects the local RBA block applied to the common
local lifted basis register. Equivalently, the circuit implements a
block-diagonal, position-controlled operation whose $\boldsymbol{\ell}$-th
block is a block-encoding of
$A_{m,\boldsymbol{\ell}}^{\mathrm{loc,red}}$. The estimates above describe
the cost of each local block-encoding. The full position-controlled
implementation may require additional control or multiplexing overhead,
depending on how the local operators vary with $\boldsymbol{\ell}$.

\section{ODE example: the Lorenz system}
\label{sec:Lorenz}

In this section, we apply the RBA algorithm to the Lorenz system and work out its explicit construction for a one-timestep update, followed by its generalization to $m$-timesteps. We verify that the RBA algorithm exactly recovers the discrete time nonlinear dynamics and demonstrate how the full monomial basis and reduced basis grow under the composition of $m$ timesteps. Furthermore, we use this system to demonstrate how the success probability of block-encoding the RBA operator depends on the number of time steps, and how scaling can be used to increase this success probability.

\subsection{Problem setup}

The Lorenz system is defined by the following system of ODEs:
\begin{equation}
\label{eq:Lorenz-system}
\begin{aligned}
    \frac{dx}{dt} &= \sigma (y-x),\\
    \frac{dy}{dt} &= x(\rho-z)-y,\\
    \frac{dz}{dt} &= xy-\beta z,
\end{aligned}
\end{equation}
where $x(t)$, $y(t)$, and $z(t)$ are the state variables and $\sigma,\rho,\beta>0$ are tunable parameters. We define the associated state vector and the initial condition as
\begin{equation*}
    \boldsymbol{u}(t)
    =
    \begin{pmatrix}
    x(t)\\
    y(t)\\
    z(t)
    \end{pmatrix}, \quad
    \boldsymbol{u}^{(0)}
    =
    \boldsymbol{u}(0)
    =
    \begin{pmatrix}
    x_0\\
    y_0\\
    z_0
    \end{pmatrix}.
\end{equation*}

The right-hand side of the Lorenz system in Eq.~\eqref{eq:Lorenz-system} is a polynomial map of degree $p=2$, due to the quadratic terms $xz$ and $xy$. Applying the explicit Euler method with timestep $h>0$ gives
\begin{equation}
\label{eq:Lorenz-map}
\begin{aligned}
    x^{(k+1)} &= \alpha x^{(k)}+\gamma y^{(k)},\\
    y^{(k+1)} &= \delta x^{(k)}+\eta y^{(k)}-h x^{(k)}z^{(k)},\\
    z^{(k+1)} &= \theta z^{(k)}+h x^{(k)}y^{(k)},
\end{aligned}
\end{equation}
where $\alpha := 1-h\sigma$, $\gamma := h\sigma$, $\delta := h\rho$, $\eta := 1-h$, and $\theta := 1-h\beta$. Thus, the one-timestep Euler update can be written as the quadratic polynomial map
\begin{equation}
\label{eq:Lorenz-one-step-map}
\boldsymbol{\Phi}_h:\mathbb{R}^3\to\mathbb{R}^3,
\qquad
\boldsymbol{\Phi}_h(\boldsymbol{u})
=
\begin{pmatrix}
    \alpha x+\gamma y\\
    \delta x+\eta y-hxz\\
    \theta z+hxy
\end{pmatrix}.
\end{equation}
The corresponding $m$-timestep composed map is denoted by $\boldsymbol{\Phi}_{h}^{\circ m}$, and evaluating this map at the initial condition gives the state after $m$ Euler timesteps:
\begin{equation}
\label{eq:Lorenz-m-step-state}
     \boldsymbol{u}^{(m)}
     =
     \boldsymbol{\Phi}_{h}^{\circ m}(\boldsymbol{u}^{(0)}).
\end{equation}
Since $\boldsymbol{\Phi}_h$ is a polynomial map of degree $2$, each component of the composed map $\boldsymbol{\Phi}_{h}^{\circ m}$ is a polynomial in the initial state variables of degree at most $2^m$.

\subsection{One-timestep RBA construction}

For a one-timestep update, $m=1$, the formal degree bound is $2$. Using the definition of the full monomial basis given in Eq. \eqref{eq:full-monomial-basis-ODE}, with $n=3$ since the Lorenz system has three state variables, the full monomial basis of total degree at most $2$ is
\begin{equation}
\label{eq:Lorenz-full basis}
\mathcal{M}_2^{\mathrm{full}}
=
\{1,x,y,z,x^2,xy,xz,y^2,yz,z^2\},
\end{equation}
with dimension
\begin{equation*}
\label{eq:Lorenz-full basis-dim}
d_2^{\mathrm{full}}
=
\left|\mathcal{M}_2^{\mathrm{full}}\right|
=
\binom{3+2}{2}
=
10.
\end{equation*}
Following the ordering in Eq. \eqref{eq:Lorenz-full basis}, the full lifted state vector evaluated at the initial condition is
\begin{equation}
\label{eq:Lorenz-full-lifted-vector-initial}
\boldsymbol{w}_1^{\mathrm{full}}(\boldsymbol{u}^{(0)})
=
\begin{pmatrix}
1\\
x_0\\
y_0\\
z_0\\
x_0^2\\
x_0y_0\\
x_0z_0\\
y_0^2\\
y_0z_0\\
z_0^2
\end{pmatrix}.
\end{equation}

Following the discrete time system given in Eq. \eqref{eq:Lorenz-map}, the one-timestep full basis RBA operator $\boldsymbol{A}_1^{\mathrm{full}}\in\mathbb{R}^{10\times 10}$ is defined as 
\begin{equation*}
\label{eq:Lorenz-full-A}
\boldsymbol{A}_1^{\mathrm{full}}
=
\begin{pmatrix}
0 & \alpha & \gamma & 0      & 0 & 0 & 0  & 0 & 0 & 0\\
0 & \delta & \eta   & 0      & 0 & 0 & -h & 0 & 0 & 0\\
0 & 0      & 0      & \theta & 0 & h & 0  & 0 & 0 & 0\\
0 & 0      & 0      & 0      & 0 & 0 & 0  & 0 & 0 & 0\\
0 & 0      & 0      & 0      & 0 & 0 & 0  & 0 & 0 & 0\\
0 & 0      & 0      & 0      & 0 & 0 & 0  & 0 & 0 & 0\\
0 & 0      & 0      & 0      & 0 & 0 & 0  & 0 & 0 & 0\\
0 & 0      & 0      & 0      & 0 & 0 & 0  & 0 & 0 & 0\\
0 & 0      & 0      & 0      & 0 & 0 & 0  & 0 & 0 & 0\\
0 & 0      & 0      & 0      & 0 & 0 & 0  & 0 & 0 & 0
\end{pmatrix}.
\end{equation*}
When applied to the full lifted vector given in Eq. \eqref{eq:Lorenz-full-lifted-vector-initial} this operator returns the state after one timestep of the Euler map, embedded in the first three entries:
\begin{equation}
\label{eq:Lorenz-full-action}
\boldsymbol{A}_1^{\mathrm{full}}
\boldsymbol{w}_1^{\mathrm{full}}(\boldsymbol{u}^{(0)})
=
\begin{pmatrix}
\boldsymbol{u}^{(1)}\\
\mathbf{0}_{d_2^{\mathrm{full}}-3}
\end{pmatrix}.
\end{equation}

This exactly represents one Euler timestep. However, the full lifted basis is not minimal. Only five distinct monomials appear in the one-timestep Euler map, as given in Eq. \eqref{eq:Lorenz-one-step-map}, namely $\{x,y,z,xz,xy\}$. Thus, we may define the reduced lifted basis
\begin{equation}
\label{eq:Lorenz-reduced basis}
\mathcal{M}_2^{\mathrm{red}}
=
\{x,y,z,xz,xy\}.
\end{equation}
The corresponding reduced lifted vector evaluated at the initial condition is
\begin{equation}
\label{eq:Lorenz-reduced-lifted-vector-initial}
\boldsymbol{w}_1^{\mathrm{red}}(\boldsymbol{u}^{(0)}\\)
=
\begin{pmatrix}
x_0\\
y_0\\
z_0\\
x_0z_0\\
x_0y_0
\end{pmatrix}.
\end{equation}
The reduced lifted dimension is $d_2^{\mathrm{red}}=5$, rather than the full basis dimension $d_2^{\mathrm{full}}=10$. With the ordering in Eq. \eqref{eq:Lorenz-reduced basis}, the reduced RBA operator is
\begin{equation*}
\label{eq:Lorenz-reduced-A}
\boldsymbol{A}_1^{\mathrm{red}}
=
\begin{pmatrix}
\alpha & \gamma & 0      & 0  & 0\\
\delta & \eta   & 0      & -h & 0\\
0      & 0      & \theta & 0  & h\\
0      & 0      & 0      & 0  & 0\\
0      & 0      & 0      & 0  & 0
\end{pmatrix}.
\end{equation*}
When applied to the reduced lifted vector defined in Eq. \eqref{eq:Lorenz-reduced-lifted-vector-initial}, this operator gives the state after one timestep of the Euler map, analogously to Eq. \eqref{eq:Lorenz-full-action}:
\begin{equation*}
\label{eq:Lorenz-reduced-action}
\boldsymbol{A}_1^{\mathrm{red}}
\boldsymbol{w}_1^{\mathrm{red}}(\boldsymbol{u}^{(0)})
=
\begin{pmatrix}
\boldsymbol{u}^{(1)}\\
\mathbf{0}_{d_2^{\mathrm{red}}-3}
\end{pmatrix}.
\end{equation*}
The reduced RBA operator gives the same one-timestep Euler update as the full basis operator while acting on a lower-dimensional lifted vector.

\subsection{Generalization to $m$-timesteps}

As shown in the previous section, the $m$-fold composition $\boldsymbol{\Phi}_{h}^{\circ m}$ for the Lorenz system is a polynomial in the initial state variables of degree at most $2^m$. Therefore, using the definition in Eq.~\eqref{eq:full-monomial-basis-ODE}, the full $m$-timestep monomial basis is given by
\begin{equation*}
\label{eq:Lorenz-m-full basis}
\mathcal{M}_{2^m}^{\mathrm{full}}
=
\left\{
x^{\nu_1}y^{\nu_2}z^{\nu_3}
:
\boldsymbol{\nu}
\in
\mathbb{N}_0^3,
\quad
|\boldsymbol{\nu}|
:=
\nu_1+\nu_2+\nu_3
\le 2^m
\right\}.
\end{equation*}
The reduced $m$-timestep basis, $\mathcal{M}_{2^m}^{\mathrm{red}}$, is defined by retaining only those monomials that appear with a nonzero coefficient in at least one component of the composed map $\boldsymbol{\Phi}_{h}^{\circ m}$. Therefore, the dimension of the reduced basis satisfies
\begin{equation*}
\label{eq:Lorenz-m-reduced-dim-bound}
d_{2^m}^{\mathrm{red}}
\le
\left|\mathcal{M}_{2^m}^{\mathrm{full}}\right|
=
\binom{3+2^m}{2^m}.
\end{equation*}

We choose an ordering of the reduced basis such that the state variables appear first,
\begin{equation*}
\label{eq:Lorenz-m-reduced basis}
    \mathcal{M}_{2^m}^{\mathrm{red}}
    =
    \{\psi_1^{\mathrm{red}},\ldots,\psi^{\mathrm{red}}_{d_{2^m}^{\mathrm{red}}}\},
    \quad
    \psi_1^{\mathrm{red}}=x,\quad
    \psi_2^{\mathrm{red}}=y,\quad
    \psi_3^{\mathrm{red}}=z,
\end{equation*}
where each $\psi_j^{\mathrm{red}}$ denotes one of the monomials in this basis. The corresponding lifted state vector evaluated at the initial condition is
\begin{equation}
\label{eq:Lorenz-m-lifted-vector}
\boldsymbol{w}_m^{\mathrm{red}}(\boldsymbol{u}^{(0)})
=
\begin{pmatrix}
\psi_1^{\mathrm{red}}(\boldsymbol{u}^{(0)})\\
\vdots\\
\psi^{\mathrm{red}}_{d_{2^m}^{\mathrm{red}}}(\boldsymbol{u}^{(0)})
\end{pmatrix}.
\end{equation}
Since the reduced basis contains all monomials appearing in the composed map, there exist scalar coefficients $c_{\alpha j}^{(m)}$ such that
\begin{equation*}
\label{eq:Lorenz-m-expansion}
\bigl(\boldsymbol{\Phi}_{h}^{\circ m}(\boldsymbol{u})\bigr)_\alpha
=
\sum_{j=1}^{d_{2^m}^{\mathrm{red}}}
c_{\alpha j}^{(m)}
\psi^{\mathrm{red}}_j(\boldsymbol{u}),
\qquad
\alpha=1,2,3.
\end{equation*}
The corresponding $m$-timestep RBA operator $\boldsymbol{A}_m^{\mathrm{red}}
\in
\mathbb{R}^{d_{2^m}^{\mathrm{red}}\times d_{2^m}^{\mathrm{red}}}$ is defined by placing these coefficients in the first three rows:
\begin{equation*}
\label{eq:Lorenz-m-A-entries}
(\boldsymbol{A}_m^{\mathrm{red}})_{\alpha j}
=
\begin{cases}
c_{\alpha j}^{(m)}, & 1\le \alpha\le 3,\\
0, & 3<\alpha\le d_{2^m}^{\mathrm{red}},
\end{cases}
\qquad
j=1,\ldots,d_{2^m}^{\mathrm{red}}.
\end{equation*}
Acting on the reduced lifted state vector defined in Eq.~\eqref{eq:Lorenz-m-lifted-vector} gives the state after $m$ timesteps of the Euler map in the first three entries:
\begin{equation*}
    \boldsymbol{A}_m^{\mathrm{red}}
    \boldsymbol{w}_m^{\mathrm{red}}(\boldsymbol{u}^{(0)})
    =
    \begin{pmatrix}
    \boldsymbol{u}^{(m)}\\
    \boldsymbol{0}_{d_{2^m}^{\mathrm{red}}-3}
    \end{pmatrix}.
\end{equation*}
Thus, for any fixed number of timesteps $m$, the RBA algorithm gives an exact linear representation of the $m$-timestep explicit Euler map. 

\subsection{Numerical Verification}

To verify numerically that the RBA reproduces the same trajectory as the corresponding explicit Euler method, we simulate the Lorenz system with parameters $\sigma = 10.0$, $\rho = 28.0$, $\beta = \frac{8.0}{3.0}$, and integrate up to final time $T=30.0$ with timestep size $h=0.001$. This corresponds to $30,000$ timesteps. In principle, one could construct an RBA operator representing the complete 30,000 timestep Euler map. However, this is not computationally feasible due to the algebraic complexity of the composed polynomial map, which grows exponentially like $2^m$. While the reduced basis is much smaller, it still grows rapidly with $m$, with the observed growth over the tested range consistent with exponential scaling as shown in Table \ref{tab:full-reduced-comparison-lorenz}. By $m=7$, the full basis already contains 366,145 monomials, while the reduced basis contains 3371 monomials. The last column reports the qubit count required to encode the reduced basis, which grows approximately linearly over the reported range.

\begin{table}[ht]
\caption{Comparison of the full and reduced lifted basis dimensions for the Lorenz system as the number of composed timesteps $m$ increases. Here, $2^m$ is the worst-case polynomial degree of the $m$-timestep composed Euler map, $d_{2^m}^{\mathrm{full}}$ is the dimension of the full monomial basis, $d_{2^m}^{\mathrm{red}}$ is the dimension of the reduced monomial basis, and $q_m^{\mathrm{red}}$ is the number of qubits required to encode the reduced basis.}
\centering
\begin{tabular}{c c r r r r}
\hline
$m$ & $2^m$ & $d_{2^m}^{\mathrm{full}}$ & $d_{2^m}^{\mathrm{red}}$ & $q_m^{\mathrm{red}}$ \\
\hline
1 & 2   & 10      & 5    & 3  \\
2 & 4   & 35      & 12   & 4  \\
3 & 8   & 165    & 28   & 5  \\
4 & 16  & 969    & 80   & 7  \\
5 & 32  & 6545   & 255  & 8  \\
6 & 64  & 47905  & 899  & 10 \\
7 & 128 & 366145 & 3371 & 12 \\
\hline
\end{tabular}
\label{tab:full-reduced-comparison-lorenz}
\end{table}

For this reason, we do not construct the RBA operator over the entire time interval at once. Instead, we compose the Euler map over a short time window of five timesteps, construct the corresponding five-timestep reduced RBA operator and monomial basis, form the lifted representation, and apply the RBA operator to advance the solution. After each time window, the lifted representation is reinitialized using the updated physical state as the new initial condition. This procedure is repeated until the final time is reached.

\begin{figure}
    \centering
    \includegraphics[width=1.0\linewidth]{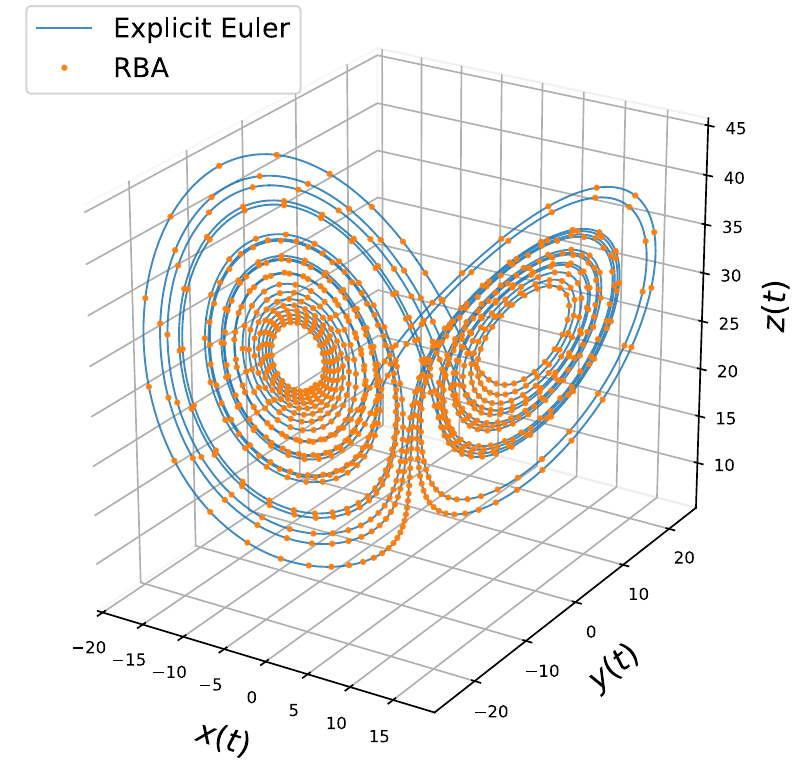}
    \caption{Numerical verification of the RBA representation for the Lorenz system, where the reference solution is obtained by the classical explicit Euler method. The two phase-space trajectories agree up to floating-point round-off.}
    \label{fig:lorenz-rba}
\end{figure}

Reinitializing after each five-timestep window prevents the lifted dimension from becoming prohibitively large while preserving the exact composed polynomial representation of the Euler updates within each window. The resulting phase-space trajectory is shown in Figure \ref{fig:lorenz-rba}, where the reference solution is obtained by the standard explicit Euler method. The two trajectories agree up to floating-point round-off, confirming that the RBA construction correctly reproduces the discrete nonlinear Euler dynamics without introducing additional error.

We also use this example to examine the post-selection probability associated with a block-encoding of the Lorenz RBA operator $\boldsymbol{A}_m^{\mathrm{red}}$. Let $\alpha_m$ denote the normalization factor used to block-encode $\boldsymbol{A}_m^{\mathrm{red}}$, such that the implemented block is $\boldsymbol{A}_m^{\mathrm{red}}/\alpha_m$, with

\begin{equation*}
\alpha_m \geq \|\boldsymbol{A}_m^{\mathrm{red}}\|_2 .
\end{equation*}
For a block initialized at timestep $k$, with $k=0$ denoting the initial condition and $k>0$ denoting reinitialization from a later state along the
trajectory, the normalized lifted state is prepared according to Eq. \eqref{eq:ODE-amplitude-encoding}. After applying the block-encoded operator and successfully post-selecting on the block-encoding ancilla, the resulting state is proportional to
\begin{equation*}
\frac{
\boldsymbol{A}_m^{\mathrm{red}}
\boldsymbol{w}_m^{\mathrm{red}}(\boldsymbol{u}^{(k)})
}{
\alpha_m
\|\boldsymbol{w}_m^{\mathrm{red}}(\boldsymbol{u}^{(k)})\|_2
}.
\end{equation*}
Therefore the post-selection success probability is
\begin{equation}
\label{eq:Lorenz-block-encoding-success}
p_{\mathrm{succ}}^{(m)}(\boldsymbol{u}^{(k)})
=
\frac{
\left\|
\boldsymbol{A}_m^{\mathrm{red}}
\boldsymbol{w}_m^{\mathrm{red}}(\boldsymbol{u}^{(k)})
\right\|_2^2
}{
\alpha_m^2
\left\|
\boldsymbol{w}_m^{\mathrm{red}}(\boldsymbol{u}^{(k)})
\right\|_2^2
}.
\end{equation}
In the numerical tests below, we take
\begin{equation*}
\alpha_m=\|\boldsymbol{A}_m^{\mathrm{red}}\|_2 .
\end{equation*}

By construction, the numerator of Eq. \eqref{eq:Lorenz-block-encoding-success} is the squared norm of the physical output obtained by applying the $m$-timestep composed Euler map. The denominator contains both the block-encoding normalization $\alpha_m$ and the squared norm of the lifted monomial vector. Hence the success probability may become small when either the matrix normalization is large or the lifted vector contains large monomial entries compared with the
physical output.

The Lorenz system demonstrates this latter effect clearly. As $m$ increases, the composed map contains monomials of increasing degree. The typical form of a monomial is $x^a y^b z^c$, whose magnitude can grow rapidly with the total degree $a+b+c$ when the physical variables are not order one. Consequently, high-degree monomial components can dominate $\|\boldsymbol{w}_m^{\mathrm{red}}(\boldsymbol{u}^{(k)})\|_2$. Since $\boldsymbol{A}_m^{\mathrm{red}}$ maps the lifted vector back to the physical three-dimensional state and sets the remaining lifted components to zero, large high-degree monomial entries can increase the denominator in Eq. \eqref{eq:Lorenz-block-encoding-success} without producing a corresponding
increase in the numerator. This explains why the block-encoding success probability decreases rapidly as the number of composed timesteps increases.

To improve the success probability, we construct an equivalent RBA representation in scaled variables. Let
\begin{equation*}
X=s_x x,\qquad Y=s_y y,\qquad Z=s_z z,
\end{equation*}
where $s_x,s_y,s_z>0$ are fixed scaling factors chosen to match the typical magnitudes of the corresponding Lorenz variables. In the scaled variables, the Euler map is obtained by transforming from $(x,y,z)$ to $(X,Y,Z)$, applying
the original Euler update, and transforming back to $(x,y,z)$. The scaled and unscaled RBA constructions describe the same explicit Euler trajectory, expressed in different coordinates. Therefore, scaling changes the conditioning of the lifted representation, and hence the post-selection probability, but it does not modify the underlying nonlinear update. Because the scaled variables are closer to order one, the high-degree monomials are typically much better balanced in magnitude. This reduces the norm of the lifted vector and therefore improves the post-selection probability.

\begin{figure}
    \centering
    \includegraphics[width=\linewidth]{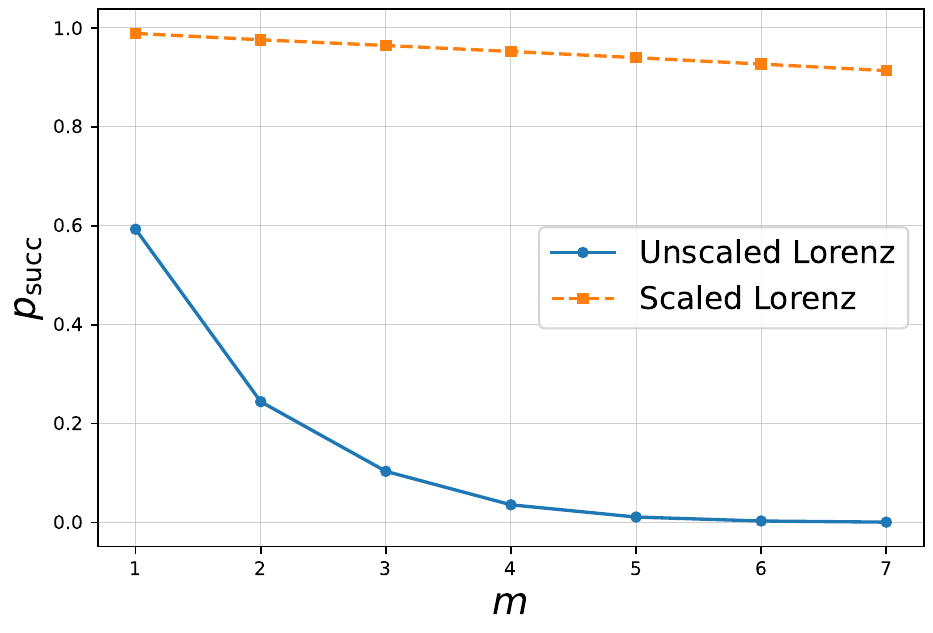}
    \caption{Block-encoding post-selection success probability for the Lorenz RBA operator as a function of the number of composed timesteps $m$. The unscaled representation uses the physical variables directly, while the scaled representation constructs the RBA map in normalized variables.}
    \label{fig:lorenz-success-probability}
\end{figure}

Figure \ref{fig:lorenz-success-probability} shows the post-selection success probabilities for the unscaled and scaled Lorenz RBA representations. In the unscaled variables, the success probability decays rapidly as the number of composed timesteps $m$ increases, indicating that the corresponding
lifted representation is poorly conditioned for block encoding. In contrast, the scaled representation substantially improves the post-selection probability, which remains close to one over the plotted range of $m$. These results demonstrate that coordinate scaling improves the
conditioning of the lifted monomial representation.

\section{PDE example: 1D viscous Burgers' equation}
\label{sec:Burgers}

In this section, we apply the RBA algorithm to the one-dimensional viscous Burgers' equation with periodic boundary conditions. We work out the construction explicitly for one timestep and then explain how the same construction extends to $m$ timesteps. Similar to the Lorenz system, we demonstrate that the RBA algorithm exactly recovers the discrete time nonlinear dynamics and demonstrate how the full monomial basis, reduced basis, and the effective stencil grow under the $m$-timestep composition.

\subsection{Problem setup}

The one-dimensional viscous Burgers' equation on the periodic domain $x\in [0,L]$ is given by
\begin{equation}
\label{eq:Burgers-system}
    \partial_t u(x,t)
    =
    -u(x,t)\partial_x u(x,t)
    +
    \nu \partial_{xx}u(x,t),
\end{equation}
where $u(x,t)$ is the scalar field of interest and $\nu\ge 0$ is the viscosity coefficient. We impose periodic boundary conditions, $u(0,t)=u(L,t)$. The nonlinear term in Eq. \eqref{eq:Burgers-system} is the quadratic self-advection term $-u(x,t)\partial_x u(x,t)$.

We discretize the domain $[0,L]$ using $N$ grid points:
\begin{equation*} 
    x_\ell=\ell\Delta x, \qquad \ell=0,\ldots,N-1, \qquad \Delta x=\frac{L}{N}. 
\end{equation*} 
Using centered finite differences, we approximate the spatial derivatives as  
\begin{equation*} 
    \partial_x u(x_\ell,t) \approx \frac{u_{\ell+1}(t)-u_{\ell-1}(t)}{2\Delta x}, 
\end{equation*} and 
\begin{equation*} 
    \partial_{xx}u(x_\ell,t) \approx \frac{u_{\ell+1}(t)-2u_\ell(t)+u_{\ell-1}(t)}{\Delta x^2}. 
\end{equation*}
For $\ell=0,\ldots,N-1$, this gives the semi-discrete ODE system
\begin{equation}
\label{eq:Burgers-semidiscrete-system}
    \frac{d u_\ell}{dt}
    =
    -u_\ell
    \frac{u_{\ell+1}-u_{\ell-1}}{2\Delta x}
    +
    \nu
    \frac{u_{\ell+1}-2u_\ell+u_{\ell-1}}{\Delta x^2},
\end{equation}
where all indices are interpreted $\mathrm{mod } \ N$.

We collect all grid values into the semi-discrete state vector
\begin{equation*}
    \label{eq:statevector-burgers}
    \boldsymbol{u}(t)
    =
    \begin{pmatrix}
    u_0(t)\\
    u_1(t)\\
    \vdots\\
    u_{N-1}(t)
    \end{pmatrix},
    \qquad
    \boldsymbol{u}^{(0)}
    =
    \boldsymbol{u}(0)
    =
    \begin{pmatrix}
    u_{\mathrm{init},x_0}\\
    u_{\mathrm{init},x_1}\\
    \vdots\\
    u_{\mathrm{init},x_{N-1}}
    \end{pmatrix}.
\end{equation*}
Then, Eq. \eqref{eq:Burgers-semidiscrete-system} can be written as a finite-dimensional polynomial ODE system, following the formulation in Eq. \eqref{eq:PDEtoODE}. Applying the explicit Euler method with timestep $h>0$ gives
\begin{equation}
\label{eq:Burgers-Euler-map}
\begin{aligned}
    u_\ell^{(k+1)}
    &=
    (1-2b)u_\ell^{(k)}
    +
    b u_{\ell-1}^{(k)}
    +
    b u_{\ell+1}^{(k)} \\
    &\quad
    +
    a u_\ell^{(k)}u_{\ell-1}^{(k)}
    -
    a u_\ell^{(k)}u_{\ell+1}^{(k)},
\end{aligned}
\end{equation}
where
\begin{equation*}
    a:=\frac{h}{2\Delta x},
    \qquad
    b:=\frac{\nu h}{\Delta x^2}.
\end{equation*}
Thus, the one-timestep Euler update in Eq. \eqref{eq:Burgers-Euler-map} can be written as the quadratic polynomial map
\begin{equation*}
    \boldsymbol{\Phi}_h:\mathbb{R}^N\to\mathbb{R}^N,
    \qquad
    \boldsymbol{u}^{(k+1)}
    =
    \boldsymbol{\Phi}_h\bigl(\boldsymbol{u}^{(k)}\bigr).
\end{equation*}
The corresponding $m$-timestep composed map is denoted by $\boldsymbol{\Phi}_{h}^{\circ m}$, and evaluating this map at the initial condition gives the state after $m$ Euler timesteps:
\begin{equation*}
     \boldsymbol{u}^{(m)}
     =
     \boldsymbol{\Phi}_{h}^{\circ m}(\boldsymbol{u}^{(0)}).
\end{equation*}
Since $\boldsymbol{\Phi}_h$ is a polynomial map of degree $2$, each component of the composed map $\boldsymbol{\Phi}_{h}^{\circ m}$ is a polynomial in the initial state variables of degree at most $2^m$.

\subsection{One-timestep RBA construction}

For a one-timestep update, $m=1$, the formal degree bound is $2$. Using the definition of the full monomial basis in Eq. \eqref{eq:PDE-global-full basis}, with $n=N$ because the spatial discretization results in an $N$-dimensional ODE system, the full monomial basis of total degree at most $2$ is
\begin{equation*}
\label{eq:Burgers-global-full basis}
    \mathcal{M}_{2}^{\mathrm{glob,full}}
    =
    \left\{
    \boldsymbol{u}^{\boldsymbol{\nu}}
    :
    \boldsymbol{\nu}\in\mathbb{N}_0^N,
    \ |\boldsymbol{\nu}|\le 2
    \right\},
    \qquad
    \boldsymbol{u}^{\boldsymbol{\nu}}
    :=
    \prod_{i=0}^{N-1} u_i^{\nu_i},
\end{equation*}
with dimension
\begin{equation*}
    d_{2}^{\mathrm{glob,full}}
    =
    \left|\mathcal{M}_{2}^{\mathrm{glob,full}}\right|
    =
    \binom{N+2}{2}.
\end{equation*}

This global basis is unnecessarily large since the one-timestep update at each grid point $\ell$ depends only on the three-point stencil
\begin{equation}
\label{eq:Burgers-one-step-stencil}
    S_1(\ell)
    =
    \{\ell-1,\ell,\ell+1\}.
\end{equation}
For $N\ge 3$, the dimension of this stencil is given by
\begin{equation*}
    s_1(\ell)= |S_1(\ell)|=3.
\end{equation*}
Since Burgers' equation has a single scalar field, the number of local scalar
variables is also
\begin{equation*}
    v_1(\ell)= s_1(\ell)=3.
\end{equation*}
We define these local variables evaluted at the initial condition as
\begin{equation*}
    x_-:=u_{\ell-1}^{(0)},
    \qquad
    x_0:=u_\ell^{(0)},
    \qquad
    x_+:=u_{\ell+1}^{(0)},
\end{equation*}
and collect them into the local stencil vector
\begin{equation*}
    \boldsymbol{x}^{(\ell)}
    =
    (x_-,x_0,x_+)^T.
\end{equation*}
The full local monomial basis of total degree at most $2$ is
\begin{equation*}
\label{eq:Burgers-local-full basis}
    \mathcal{M}_{2}^{\mathrm{loc,full}}(\ell)
    =
    \{
    1,
    x_-,
    x_0,
    x_+,
    x_-^2,
    x_-x_0,
    x_-x_+,
    x_0^2,
    x_0x_+,
    x_+^2
    \},
\end{equation*}
with dimension
\begin{equation*}
    d_{2}^{\mathrm{loc,full}}(\ell)
    =
    \left|
    \mathcal{M}_{2}^{\mathrm{loc,full}}(\ell)
    \right|
    =
    \binom{3+2}{2}
    =
    10.
\end{equation*}

Using Eq. \eqref{eq:Burgers-Euler-map}, we can define for each grid point $\ell$, a one-timestep local update polynomial
\begin{equation*}
\label{eq:Burgers-one-step-local-poly}
    g_{\ell}^{(1)}(x_-,x_0,x_+)
    =
    (1-2b)x_0
    +
    b x_-
    +
    b x_+
    +
    a x_0x_-
    -
    a x_0x_+,
\end{equation*}
such that
\begin{equation*}
    u_\ell^{(1)}
    =
    g_{\ell}^{(1)}(x_-,x_0,x_+).
\end{equation*}
Only five distinct monomials appear in the expression that defines $g_{\ell}^{(1)}$, therefore, the reduced local monomial basis is
\begin{equation*}
\label{eq:Burgers-one-step-reduced basis}
    \mathcal{M}_{2}^{\mathrm{loc,red}}(\ell)
    =
    \{
    x_0,
    x_-,
    x_+,
    x_0x_-,
    x_0x_+
    \},
\end{equation*}
with dimension
\begin{equation*}
    d_{2}^{\mathrm{loc,red}}(\ell)
    =
    \left|
    \mathcal{M}_{2}^{\mathrm{loc,red}}(\ell)
    \right|
    =
    5.
\end{equation*}

The ordering of this monomial basis is chosen such that the first monomial is the central
state variable $x_0=u_\ell^{(0)}$. The corresponding reduced local lifted vector evaluated on the initial stencil is
\begin{equation}
\label{eq:Burgers-one-step-lifted-vector}
    \boldsymbol{w}_{1,\ell}^{\mathrm{loc,red}}
    \bigl(\boldsymbol{x}^{(\ell)}\bigr)
    =
    \begin{pmatrix}
        x_0\\
        x_-\\
        x_+\\
        x_0x_-\\
        x_0x_+
    \end{pmatrix}.
\end{equation}

The one-timestep local RBA operator
$\boldsymbol{A}_{1,\ell}^{\mathrm{loc,red}}\in\mathbb{R}^{5\times 5}$ is
defined as
\begin{equation*}
\label{eq:Burgers-one-step-local-A}
    \boldsymbol{A}_{1,\ell}^{\mathrm{loc,red}}
    =
    \begin{pmatrix}
        1-2b & b & b & a & -a\\
        0 & 0 & 0 & 0 & 0\\
        0 & 0 & 0 & 0 & 0\\
        0 & 0 & 0 & 0 & 0\\
        0 & 0 & 0 & 0 & 0
    \end{pmatrix}.
\end{equation*}
Although we have written this operator with an $\ell$ subscript to emphasize that it acts on the local stencil centered at grid point $\ell$, its matrix entries do not depend on $\ell$. This is a consequence of the periodic boundary conditions and the translation-invariance of the finite-difference stencil. Hence, $\boldsymbol{A}_{1,\ell}^{\mathrm{loc,red}} = \boldsymbol{A}_{1}^{\mathrm{loc,red}}$ for all $\ell$ . The grid-point dependence enters only through the local stencil variables
used to form the lifted vector in Eq. \eqref{eq:Burgers-one-step-lifted-vector}. Acting on the reduced local lifted vector in
Eq. \eqref{eq:Burgers-one-step-lifted-vector} gives
\begin{equation*}
\label{eq:Burgers-one-step-local-action}
    \boldsymbol{A}_{1}^{\mathrm{loc,red}}
    \boldsymbol{w}_{1,\ell}^{\mathrm{loc,red}}
    \bigl(\boldsymbol{x}^{(\ell)}\bigr)
    =
    \begin{pmatrix}
        u_\ell^{(1)}\\
        \boldsymbol{0}_{d_{2}^{\mathrm{loc,red}}-1}
    \end{pmatrix}.
\end{equation*}
Thus, the nonlinear one-timestep update is represented exactly as a linear map acting on the reduced local lifted monomial vector.

Since the same local RBA matrix is used at every grid point, the quantum algorithm implementation does not require a nontrivial position-dependent
control on the local RBA operator. To write the corresponding global RBA operator, we order the global lifted
vector by concatenating the reduced local lifted vectors over all grid points:
\begin{equation*}
\boldsymbol W_1^{\mathrm{loc,red}}(\boldsymbol u^{(0)})
=
\begin{pmatrix}
\boldsymbol w_{1,0}^{\mathrm{loc,red}}
\bigl(\boldsymbol x^{(0)}\bigr)\\
\vdots\\
\boldsymbol w_{1,N-1}^{\mathrm{loc,red}}
\bigl(\boldsymbol x^{(N-1)}\bigr)
\end{pmatrix}.
\end{equation*}
With this ordering, the global RBA operator is block diagonal, with one
identical local RBA matrix acting on each grid-point block:
\begin{equation}
\label{eq:global-A-burgers}
\mathcal A_1^{\mathrm{red}}
=
I_N\otimes \boldsymbol A_1^{\mathrm{loc,red}}.
\end{equation}
Equivalently,
\begin{equation*}
\mathcal A_1^{\mathrm{red}}
\boldsymbol W_1^{\mathrm{loc,red}}(\boldsymbol u^{(0)})
=
\begin{pmatrix}
u_0^{(1)}\\
\boldsymbol 0\\
\vdots\\
u_{N-1}^{(1)}\\
\boldsymbol 0
\end{pmatrix},
\end{equation*}
where, in each local block, the first entry is the updated physical value and
the remaining entries are zero. The physical solution vector
$\boldsymbol u^{(1)}$ is obtained by selecting the first component of each
local output block.
\subsection{Generalization to $m$-timesteps}

The one-timestep Burgers map in Eq. \eqref{eq:Burgers-Euler-map} is a local quadratic polynomial map. At grid point $\ell$, the one-timestep update depends only on the stencil defined in Eq. \eqref{eq:Burgers-one-step-stencil}. After composing the map with itself $m$ times, the dependence region expands because each neighboring value appearing after one timestep depends on its own neighboring
values at the previous timestep. Hence, after $m$ timesteps, the value $u_\ell^{(m)}$ depends only on variables in the effective stencil
\begin{equation*}
\label{eq:Burgers-m-step-stencil}
    S_m(\ell)
    =
    \{(\ell+r)\bmod N : r=-m,\ldots,m\},
\end{equation*}
whose dimension is given by
\begin{equation*}
    s_m(\ell)
    =
    |S_m(\ell)|
    =
    \min(N,2m+1).
\end{equation*}

In the non-wrapped regime $2m+1\le N$, this reduces to
\begin{equation*}
    s_m(\ell)=2m+1.
\end{equation*}
Since Burgers' equation has one scalar field, the number of local scalar
variables is
\begin{equation*}
    v_m(\ell)=s_m(\ell)=2m+1.
\end{equation*}
In this regime, we define the local variables evaluated at the initial condition as
\begin{equation*}
    x_r:=u_{\ell+r}^{(0)},
    \qquad
    r=-m,\ldots,m,
\end{equation*}
and collect them into the local vector
\begin{equation*}
    \boldsymbol{x}^{(\ell)}
    =
    (x_{-m},x_{-m+1},\ldots,x_0,\ldots,x_{m-1},x_m)^T
    \in\mathbb{R}^{2m+1}.
\end{equation*}
If the stencil wraps around the periodic domain, repeated periodic indices are
identified and $\boldsymbol{x}^{(\ell)}$ contains the $s_m(\ell)$ distinct
variables in $S_m(\ell)$.

Since the one-timestep map is quadratic, the formal degree bound after $m$
timesteps is $2^m$. Therefore, there exists a local polynomial
\begin{equation*}
    g_{\ell}^{(m)}
    :
    \mathbb{R}^{v_m(\ell)}
    \to
    \mathbb{R},
    \qquad
    \deg g_{\ell}^{(m)}\le 2^m,
\end{equation*}
such that
\begin{equation*}
\label{eq:Burgers-m-step-local-map}
    u_\ell^{(m)}
    =
    g_{\ell}^{(m)}
    \bigl(\boldsymbol{x}^{(\ell)}\bigr).
\end{equation*}
Equivalently, $g_{\ell}^{(m)}$
is the $\ell$-th component of the composed Euler map $\boldsymbol{\Phi}_{h}^{\circ m}$, restricted to the variables in the effective
stencil $S_m(\ell)$. The full local monomial basis for the $m$-timestep update is
\begin{equation*}
\label{eq:Burgers-m-full-local-basis}
    \mathcal{M}_{2^m}^{\mathrm{loc,full}}(\ell)
    =
    \left\{
    \bigl(\boldsymbol{x}^{(\ell)}\bigr)^{\boldsymbol{\nu}}
    :
    \boldsymbol{\nu}\in\mathbb{N}_0^{v_m(\ell)},
    \ |\boldsymbol{\nu}|\le 2^m
    \right\},
\end{equation*}
with dimension
\begin{equation*}
\label{eq:Burgers-m-full-local-dim}
    d_{2^m}^{\mathrm{loc,full}}(\ell)
    =
    \left|
    \mathcal{M}_{2^m}^{\mathrm{loc,full}}(\ell)
    \right|
    =
    \binom{v_m(\ell)+2^m}{2^m}.
\end{equation*}
In the non-wrapped regime $2m+1\le N$, this becomes
\begin{equation*}
    d_{2^m}^{\mathrm{loc,full}}(\ell)
    =
    \binom{2m+1+2^m}{2^m}.
\end{equation*}

As for the one-timestep construction, the full basis is not minimal. We therefore define the reduced local basis by retaining only the monomials that appear with nonzero coefficient in the local polynomials $g_{\ell}^{(m)}$ such that $\mathcal{M}_{2^m}^{\mathrm{loc,red}}(\ell)
    \subseteq
    \mathcal{M}_{2^m}^{\mathrm{loc,full}}(\ell)$.
If we denote its dimension by
\begin{equation*}
    d_{2^m}^{\mathrm{loc,red}}(\ell)
    =
    \left|
    \mathcal{M}_{2^m}^{\mathrm{loc,red}}(\ell)
    \right|,
\end{equation*}
then
\begin{equation*}
\label{eq:Burgers-m-reduced-local-dim-bound}
    d_{2^m}^{\mathrm{loc,red}}(\ell)
    \le
    d_{2^m}^{\mathrm{loc,full}}(\ell)
    =
    \binom{v_m(\ell)+2^m}{2^m}.
\end{equation*}

After choosing an ordering of the reduced local basis
\begin{equation*}
    \mathcal{M}_{2^m}^{\mathrm{loc,red}}(\ell)
    =
    \left\{
    \psi_{1,\ell}^{\mathrm{loc}},
    \ldots,
    \psi_{d_{2^m}^{\mathrm{loc,red}}(\ell),\ell}^{\mathrm{loc}}
    \right\},
\end{equation*}
the corresponding reduced local lifted vector is
\begin{equation}
\label{eq:Burgers-m-local-lifted-vector}
    \boldsymbol{w}_{m,\ell}^{\mathrm{loc,red}}
    \bigl(\boldsymbol{x}^{(\ell)}\bigr)
    =
    \begin{pmatrix}
    \psi_{1,\ell}^{\mathrm{loc}}
    \bigl(\boldsymbol{x}^{(\ell)}\bigr)\\
    \vdots\\
    \psi_{d_{2^m}^{\mathrm{loc,red}}(\ell),\ell}^{\mathrm{loc}}
    \bigl(\boldsymbol{x}^{(\ell)}\bigr)
    \end{pmatrix}.
\end{equation}
Since the reduced basis contains all monomials appearing in
$g_{\ell}^{(m)}$, there exist coefficients $c_j^{(m,\ell)}$ such that
\begin{equation}
\label{eq:Burgers-m-local-expansion}
    g_{\ell}^{(m)}
    \bigl(\boldsymbol{x}^{(\ell)}\bigr)
    =
    \sum_{j=1}^{d_{2^m}^{\mathrm{loc,red}}(\ell)}
    c_j^{(m,\ell)}
    \psi_{j,\ell}^{\mathrm{loc}}
    \bigl(\boldsymbol{x}^{(\ell)}\bigr).
\end{equation}

The corresponding $m$-timestep local RBA operator $\boldsymbol{A}_{m,\ell}^{\mathrm{loc,red}}
    \in
    \mathbb{R}^{d_{2^m}^{\mathrm{loc,red}}(\ell)
    \times
    d_{2^m}^{\mathrm{loc,red}}(\ell)}$ is defined by placing the coefficients in the first row and setting all remaining rows to zero:
\begin{equation*}
\label{eq:Burgers-m-local-A-entries}
    \bigl(\boldsymbol{A}_{m,\ell}^{\mathrm{loc,red}}\bigr)_{ij}
    =
    \begin{cases}
    c_j^{(m,\ell)}, & i=1,\\
    0, & 2\le i\le d_{2^m}^{\mathrm{loc,red}}(\ell),
    \end{cases}
\end{equation*}
where $j=1,\ldots,d_{2^m}^{\mathrm{loc,red}}(\ell)$. Acting on the reduced local lifted vector defined in Eq. \eqref{eq:Burgers-m-local-lifted-vector} gives
\begin{equation*}
\label{eq:Burgers-m-local-action}
    \boldsymbol{A}_{m,\ell}^{\mathrm{loc,red}}
    \boldsymbol{w}_{m,\ell}^{\mathrm{loc,red}}
    \bigl(\boldsymbol{x}^{(\ell)}\bigr)
    =
    \begin{pmatrix}
    g_{\ell}^{(m)}
    \bigl(\boldsymbol{x}^{(\ell)}\bigr)\\
    \boldsymbol{0}_{d_{2^m}^{\mathrm{loc,red}}(\ell)-1}
    \end{pmatrix}
    =
    \begin{pmatrix}
    u_\ell^{(m)}\\
    \boldsymbol{0}_{d_{2^m}^{\mathrm{loc,red}}(\ell)-1}
    \end{pmatrix}.
\end{equation*}
Thus, for any fixed number of timesteps $m$, the local RBA operator gives an exact linear representation of the $m$-timestep explicit Euler update at grid point $\ell$.

The global $m$-timestep RBA operator is defined analogously to the one-timestep case. Since the periodic finite-difference stencil is translation-invariant, the local composed RBA matrix is independent of the grid point $\ell$, so $\boldsymbol A_{m,\ell}^{\mathrm{loc,red}}
= \boldsymbol A_m^{\mathrm{loc,red}}$ for all $\ell$. With the global lifted vector ordered as a concatenation of local lifted vectors, the global operator is defined as 
\begin{equation}
\label{eq:global-A-burgers-m-step}
\mathcal A_m^{\mathrm{red}}
=
I_N\otimes \boldsymbol A_m^{\mathrm{loc,red}} .
\end{equation}

\subsection{Numerical verification}
\label{subsec:Burgers-numerical-verification}

To verify that the RBA representation gives the same discrete solution as the explicit Euler finite-difference scheme for the viscous Burgers equation, we consider the periodic domain $[0,2\pi]$, with viscosity $\nu=0.03$, and initial condition $u(x,0)=\sin(x)$. The domain is discretized using $N=256$ uniformly spaced grid points, giving $\Delta x=2\pi/256$. Time integration is performed with the explicit Euler method using timestep size $h=0.0025$, and the solution is advanced to $T=0.4$, corresponding to $T/h=160$ Euler timesteps. 

For the RBA computation, the reduced local basis and corresponding RBA operator are constructed for a two-timestep Euler composition, i.e. $m=2$. The RBA trajectory is then obtained by repeatedly re-lifting the current state and applying the two-step RBA operator over successive time windows. The choice $m=2$ is made to keep the composition window small, such that the preprocessing cost remains low for this verification example.

Table \ref{tab:burgers_basis_growth} quantifies the growth of the monomial basis for the Burgers system. The reduced local basis remains much smaller than the full local monomial basis, but it still grows rapidly with $m$. For $m=2$, the effective stencil contains $s_2=5$ grid points. The full local basis has $126$ monomials, whereas the reduced basis contains only $34$ monomials, requiring $q_2^{\rm loc}=6$ qubits to encode the local lifted basis. Together with the $q_{\rm grid}=8$ qubits required for the position register for the $N=256$ grid points, this gives $q_2^{\rm PDE}=14$ qubits before block-encoding ancillas. Increasing the window to $m=3$ already raises the reduced basis to $700$ monomials and the total register size to $18$ qubits, while $m=4$ gives $63{,}424$ reduced monomials and $24$ total register qubits.

\begin{figure}
    \centering
    \includegraphics[width=\linewidth]{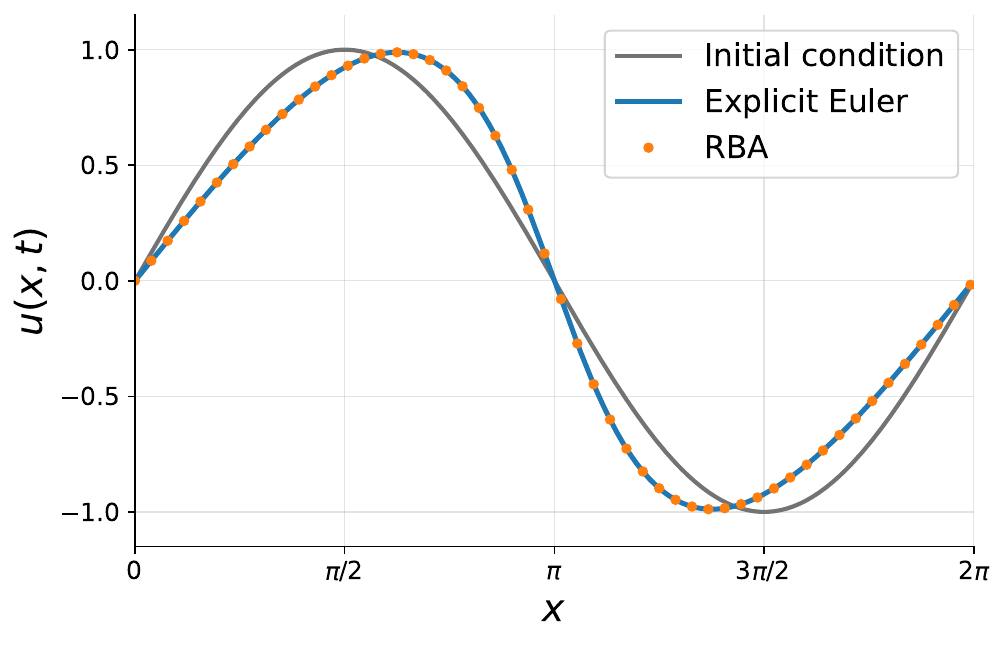}
    \caption{
    Numerical verification of the local RBA construction for the one-dimensional viscous Burgers equation. The solution is evolved to $T=0.4$ from $u(x,0)=\sin x$ using $N=256$, $\nu=0.03$, and $h=0.0025$. The RBA solution is obtained by applying the $m=2$ local lifted representation over successive two-timestep windows, with the lifted variables reconstructed from the current state at the beginning of each window. The reference solution is given by the standard finite-difference Euler implementation. The two profiles agree up to floating-point round-off.
    }
    \label{fig:burgers-rba-Euler}
\end{figure}

Using the two-timestep RBA construction, the numerical result is shown in Figure \ref{fig:burgers-rba-Euler}. The reference trajectory is obtained from the standard finite-difference Euler update. At the final time, the two solutions are visually indistinguishable and their agreement up to floating-point round-off confirms that the local RBA construction reproduces the chosen fully discrete Burgers dynamics exactly over each reinitialized time window.

\begin{table}
\caption{Growth of the local full and reduced monomial bases for the one-dimensional viscous Burgers equation. Here $m$ is the number of composed Euler steps, $2^m$ is the maximum polynomial degree of the composed map, $s_m$ is the size of the effective $m$-step stencil, $d^{\mathrm{loc,full}}_{2^m}$ is the full local monomial-basis dimension, and $d^{\mathrm{loc,red}}_{2^m}$ is the reduced local basis dimension. The last three columns give the position-register qubit count $q_{\mathrm{grid}}$, the local lifted basis qubit count $q^{\mathrm{loc}}_m$, and the total PDE register size $q^{\mathrm{PDE}}_m=q_{\mathrm{grid}}+q^{\mathrm{loc}}_m$.}
\label{tab:burgers_basis_growth}
\centering
\begin{tabular}{c c c c c c c c}
\hline
$m$ & $2^m$ & $s_m$ & $d^{\mathrm{loc,full}}_{2^m}$ 
& $d^{\mathrm{loc,red}}_{2^m}$ 
& $q_{\mathrm{grid}}$ & $q^{\mathrm{loc}}_m$ & $q^{\mathrm{PDE}}_m$ \\
\hline
1 & 2  & 3 & 10      & 5     & 8 & 3  & 11 \\
2 & 4  & 5 & 126     & 34    & 8 & 6  & 14 \\
3 & 8  & 7 & 6435    & 700   & 8 & 10 & 18 \\
4 & 16 & 9 & 2042975 & 63424 & 8 & 16 & 24 \\
\hline
\end{tabular}
\end{table}

\section{Conclusion}
\label{sec:conclusion}

In this work, we have introduced a reduced basis algorithm for simulating polynomial nonlinear ODEs and PDEs on quantum computers. The central idea is to avoid applying nonlinear transformations directly to quantum amplitudes. Instead, after time discretization, the nonlinear update map is composed over a fixed window of $m$ timesteps and represented as a polynomial in the initial state variables. By retaining only the monomials that actually appear in this composed map, the nonlinear discrete evolution can be written as a linear operator acting on a reduced lifted monomial basis.

This construction differs from mean-field and Carleman-type approaches in an important way. At the level of the chosen time discretization scheme, the RBA introduces no additional nonlinear approximation, truncation, or mean-field error. The only numerical error is the error already associated with the underlying time discretization scheme. The price for this exactness is shifted into a classical preprocessing step required to construct the composed polynomial map, identify its reduced basis, and compute the corresponding RBA operator.

For polynomial ODE systems, the lifted register size is governed by the number of monomials appearing in the $m$-timestep map. In the worst-case full basis setting, the number of qubits grows at most linearly in $m$ for fixed state dimension and polynomial degree, although the lifted basis dimension itself may grow exponentially. For PDEs, the locality of spatial discretizations makes it possible to construct reduced bases and  RBA operators associated with local effective stencils. This preserves logarithmic dependence on the global number of grid points, while the nonlinear overhead is controlled by the growth of the local stencil and the reduced local monomial basis.

The numerical examples confirm the exactness of the construction for the nonlinear discrete dynamics. For the Lorenz system, the RBA reproduces the explicit Euler trajectory up to floating-point round-off when applied over repeated short time windows. For the one-dimensional viscous Burgers equation, the local RBA construction likewise reproduces the finite-difference Euler solution over successive reinitialized time windows. These examples also demonstrate the main practical limitation of the method, which is that although the reduced basis can be much smaller than the full monomial basis, it can still grow rapidly with the number of composed timesteps. Consequently, the method is naturally suited to short or moderate time windows, with longer simulations obtained through repeated reinitialization.

A central direction for future work is to identify problem classes for which the RBA algorithm has a realistic path toward quantum advantage. The most promising options are likely to be those in which the reduced basis dimension grows slowly with the number of composed timesteps, the polynomial degree, or the spatial discretization size. Such behavior may occur for systems with strong locality, repeated stencil structure, translation invariance, sparsity, or other algebraic constraints that limit the number of distinct monomials generated under composition.

Another important direction is to determine when the classical RBA preprocessing can be reused across a family of simulations. For example, in
PDE applications, many simulations may share the same local differential operator and discretization scheme, while differing only in
initial data, forcing, boundary conditions, or parameter values. In such cases, one may not need to construct a new RBA polynomial map for every simulation. Instead, a single precomputed family of local RBA operators could potentially be reused across many problem instances, with only the state preparation or boundary treatment changing. This would make the preprocessing cost lower and could substantially improve the practical relevance of the algorithm.

\begin{acknowledgments}
We gratefully acknowledge support from the joint research project Quantum
Computational Fluid Dynamics by Fujitsu Limited and Delft University of
Technology, co-funded by the Netherlands Enterprise Agency under project
number PPS23-3-03596728.
\end{acknowledgments}

\nocite{*}
\bibliography{aipsamp}

\end{document}